\documentclass[onefignum,onetabnum]{siamart190516}

\AtBeginDocument{%
  \providecommand\BibTeX{{%
    \normalfont B\kern-0.5em{\scshape i\kern-0.25em b}\kern-0.8em\TeX}}}

\usepackage[utf8]{inputenc}
\usepackage{amsfonts}
\usepackage{amsmath}
\usepackage{tikz}
\usepackage{algorithm}
\usepackage{algorithmic}
\usepackage[export]{adjustbox}
\usepackage{svg}
\usepackage{pdfpages}
\usepackage{tabularx}
\usepackage{wrapfig}
\usepackage{array}

\newcommand{\Sum}{\sum\limits}
\newcommand{\norm}[1]{\left\lVert#1\right\rVert}

\newcommand{\midd}{\; \mid \;}
\newcommand{\set}[1]{\left\{#1\right\}}
\newcommand{\RR}{\mathbb{R}}

\newcommand{\phisidebarycenter}[1]{\int_{x\in C_{j, i_k}}x d\mu(x)}
\newcommand{\phisidearea}[1]{\int_{x\in C_{j, i_k}}d\mu(x)}
\newcommand{\Lagx}{Pow^\phi(x_i)}
\newcommand{\Lagy}{Pow^\psi(y_i)}

\title{Symmetrized semi-discrete optimal transport}

\author{Agathe Herrou\thanks{CNRS, UCBL, LIRIS (\email{agathe.herrou@liris.cnrs.fr})},
\and Bruno Lévy\thanks{INRIA Nancy Grand Est (\email{bruno.levy@inria.fr})},
\and Vincent Nivoliers\thanks{UCBL, CNRS, LIRIS (\email{vincent.nivoliers@liris.cnrs.fr})},
\and Nicolas Bonneel\thanks{CNRS, UCBL, LIRIS (\email{nicolas.bonneel@liris.cnrs.fr})},
\and Julie Digne\thanks{CNRS, UCBL, LIRIS (\email{julie.digne@liris.cnrs.fr})}}

\graphicspath{{fig/}}

\begin{document}

\maketitle

\begin{abstract}
  Interpolating between measures supported by polygonal or polyhedral domains is a problem that has been recently
  addressed by the semi-discrete optimal transport framework. Within this framework, one of the domains is discretized with a set of samples, while the other one remains continuous.
  In this paper we present a method to introduce some symmetry into the solution using coupled power diagrams. This symmetry is key to capturing the discontinuities of the transport map reflected in the geometry of the power cells.
  We design our method as a fixed-point algorithm alternating between computations of semi-discrete transport maps and recentering of the sites. 
  The resulting objects are coupled power diagrams with identical geometry, allowing us to approximate displacement interpolation through linear interpolation of the meshes' vertices.
  Through these coupled power diagrams, we have a natural way of jointly sampling measures.

\end{abstract}

\section{Introduction}
Interpolating between mathematical distributions has a 
wide range of applications from histogram
interpolation to mesh morphing. Different notions of interpolation have been 
proposed, the simplest certainly being linear interpolation. However, linear interpolation is not adapted to all applications. 
In particular, when objects differ by a translation, it does not interpolate the translation itself but instead blends the two objects.
To avoid such a blending and get a real translation, researchers have turned to 
approaches such as displacement interpolation that makes use of the optimal transport framework. 
This interpolation gives more natural results for applications that require  
morphing between shapes. The downside is that optimal transport yields an optimization problem that can be very costly to solve.

Recently, semi-discrete optimal transport has produced impressive displacement interpolation results between measures supported on 
polygonal and polyhedral domains with finely discretized meshes, using efficient geometric computations based on power cells.
However, when discretization is coarse, the ``tearing'' of geometrical shapes during the displacement interpolation process  
results in visible discretization artifacts (see Fig.~\ref{fig:shape9}, second and third rows).

In this article, we present an alternating algorithm utilizing semi-discrete optimal 
transport that computes accurate and discontinuities-capturing interpolations between meshes (Fig.~\ref{fig:shape9}, top row).
To achieve this goal, we compute a semi-discrete optimal transport map from each of 
these meshes to a discrete sampling of the other, under the constraint that the location 
of each sample corresponds to the barycenter of a power cell related to the other distribution. We observe that upon convergence of our alternating minimization, 
the geometry of the power cells captures well the geometry of the transport map discontinuities -- and hence, well approximates tearing. 
A displacement interpolation can then be computed between the two resulting coupled transport maps, 
approximating the initial meshes.

\begin{figure}[!h]
    \centering
    \begin{tabular}{|c|c|c|}
      \hline
      &Input  &  Interpolation \\ 
      & measures & \\ 
      \hline   
      \rotatebox[origin=l]{90}{\scriptsize{Our interpolation} } &
      \includesvg[width=4cm]{interpolation/2d/shape9/shape9-scene.svg} &

      \includesvg[height=3cm]{interpolation/2d/shape9/shape9-step0.svg}
      \includesvg[height=3cm]{interpolation/2d/shape9/shape9-step1.svg}
      \includesvg[height=3cm]{interpolation/2d/shape9/shape9-step2.svg} \\ 
      \hline
      
     \rotatebox[origin=l]{90}{\scriptsize{semi-discrete~\cite{bl-semi-discrete}}}
      \rotatebox[origin=l]{90}{\scriptsize{$\mu$ continuous}} \rotatebox[origin=l]{90}{\scriptsize{$\nu$ sampled}} &
      \includesvg[width=4cm]{interpolation/2d/shape9/shape9-nonsym.svg} &
      
      \includesvg[height=3cm]{interpolation/2d/shape9/shape9-ns-step2.svg}
      \includesvg[height=3cm]{interpolation/2d/shape9/shape9-ns-step1.svg}
      \includesvg[height=3cm]{interpolation/2d/shape9/shape9-ns-step0.svg}\\
      \hline
      
     \rotatebox[origin=l]{90}{\scriptsize{semi-discrete~\cite{bl-semi-discrete}}}
      \rotatebox[origin=l]{90}{\scriptsize{$\mu$ sampled}} \rotatebox[origin=l]{90}{\scriptsize{$\nu$ continuous}} &
      \includesvg[width=4cm]{interpolation/2d/shape9/shape9-nonsym-rev.svg}&
      
      \includesvg[height=3cm]{interpolation/2d/shape9/shape9-ns-rev-step0.svg}
      \includesvg[height=3cm]{interpolation/2d/shape9/shape9-ns-rev-step1.svg}
      \includesvg[height=3cm]{interpolation/2d/shape9/shape9-ns-rev-step2.svg}\\
      \hline
    \end{tabular}
    \caption{Interpolation between two horizontally aligned disk and two vertically aligned disks.
      Notice how the cells align with the discontinuities on both measures with our method,
      while they only align on the target measure with Levy's method, resulting in artifacts
      in the interpolation.}\label{fig:shape9}
\end{figure}

\subsection{Previous work}
The field of optimal transport has seen numerous developments in the last decades,
as reviewed in the books of Santambrogio \cite{santambrogio2015optimal}
and Peyré and Cuturi~\cite{peyre2017computational} with several applications,
notably in Computer Graphics and Computational Geometry.
Numerical methods for discrete and semi-discrete optimal transport have been
reviewed by Mérigot and Thibert~\cite{DBLP:journals/corr/abs-2003-00855}.

In particular, semi-discrete optimal transport has been appreciated for its 
ability to bridge discrete and continuous measures, as well as for its relative ease of
computation. Semi-discrete optimal transport consists in finding an optimal transport map describing how to assign regions of a domain endowed with a continuous measure to Dirac masses -- e.g., describing how a density of population would share a set of bakeries with limited supplies. 
Aurenhammer, Hoffmann and Aronov~\cite{aha} have shown that in a semi-discrete
transport map, the pre-images of the Dirac masses correspond to a power diagram
(i.e., an extension of a Voronoi diagram with weights controlling the cells area)
with well-chosen weights, and that these weights could be expressed as the
minimum of a certain convex function.
This is also a consequence of Brenier's polar factorization theorem considered
in this particular setting~\cite{brenier-polar-factorization}.
Efficient methods have then been proposed for minimizing this function, notably a
multilevel method proposed by Mérigot~\cite{merigot-multiscale}, 
a numerical quasi-newton method by Lévy~\cite{bl-semi-discrete}
and a Newton method with proven convergence by Kitagawa et al.~\cite{kmt}.

Semi-discrete optimal transport has a wide range of applications including the computation of samplings of density functions with blue noise properties~\cite{blue-noise}, fluid simulation~\cite{gallouet2018lagrangian},  the reconstruction
of caustics for designing optical components~\cite{light-in-power}, the reconstruction
of the early stages of the universe~\cite{early-universe} from observations,
reconstructing and improving meshes~\cite{meshrec}, or the enrichment of
generated textures~\cite{enriching-textures}. 
One last application that has motivated this article is an 
interpolation framework between meshes proposed by Lévy~\cite{bl-semi-discrete}. Optimal transport is useful in this 
context where typical feature-point based interpolation would not be appropriate, as features point extraction is error-prone on simple shapes without any strongly identified local feature, such as the disks in Figure~\ref{fig:shape9}.

Mesh interpolation based on optimal transport takes advantage of the easy
computation of semi-discrete optimal transport maps by approximating one of the
meshes by a discrete set of samples and computing the optimal transport between
the other mesh and the samples.
When interpolating between two shapes using semi-discrete optimal transport,
shape A is considered to be continuous, and shape B is discretized by a set of points.
The continuous transport map from A to B is approximated by a linear interpolation between
a triangulation of the samples over B and the Delaunay triangulation that is dual to the
power diagram that represents the transport map.
Depending on the shapes to be interpolated, there can also be topology changes
such as splitting and merging of connected components. If A is continuous and B is discretized by a well
distributed set of points, then splittings are well represented, but merges can
result in jagged boundaries (Fig.~\ref{fig:shape9} second row). One can also compute the optimal 
transport from a continuous B to a sampled A. This results in the reverse effect: merges are well represented and splitting result in jagged boundaries 
(Fig.~\ref{fig:shape9} third row). Our goal is to design a method that mutually
samples A and B in such a way that both splitting and mergings are well
represented (Fig.~\ref{fig:shape9} first row). In fact, optimal transport maps
are known to be continuous when transporting any shape to a convex
shape~\cite{caffarelli1992regularity}, while transporting a convex shape to a
connected shape whose boundary has regions of sufficiently high negative total
curvature necessarily results in discontinuous transport
maps~\cite{chodosh2015discontinuity}. Finding the locus of these discontinuities
in the general case is a difficult open problem. 

Our algorithm draws its inspiration from alternating algorithms such as the one used for Centroidal Voronoi Diagram computation.
In this special type of diagram, the seed associated to each Voronoi cell is at the barycenter of the cell.
Centroidal Voronoi Tessellations were introduced by Du et al.~\cite{cvt} and studied from an optimization
point of view in Liu et al.~\cite{cvt-smoothness}. They
can be computed by using Lloyd's algorithm~\cite{lloyd}, in which one alternatively 
computes the Voronoi tessellation of a set of sites and relocates these sites to their 
cell's barycenters.
Centroidal power diagrams, studied by Xin et al.~\cite{cccpd}, extend Centroidal Voronoi Diagrams by using a power diagram with capacity constraints for the cells instead of Voronoi cells. 
Several methods have been
proposed to compute such objects \cite{Balzer2009}, including formulations through 
optimal transport \cite{blue-noise}. In either cases, a diagram computation step 
(fulfilling the capacity constraints through optimal transport or other means) is 
alternated with a recentering step. 

The algorithm presented in this article uses a similar template, and can be thought of
as both an extension of Lloyd's algorithm and the capacity-constrained centroidal power
diagram between two coupled measures instead of one, and as a symmetrization of Lévy's
mesh interpolation algorithm~\cite{bl-semi-discrete}.

\subsection{Semi-discrete optimal transport}\label{subsec:tosd}

We consider a probability measure $\mu$ supported on a subset of $\mathbb{R}^d$,  
and $N$ discrete samples $X = \set{x_i \midd 1 \le i \le N} \subset \mathbb{R}^d$,
equipped with a probability measure $\sum_{i=1}^N p_i \delta_{x_i}$ such that
$\sum_{i=1}^N p_i = \int_{x \in \mathbb{R}^d} d\mu(x)$. 
The semi-discrete optimal transport problem consists in finding a transport map 
$T: \RR^d \rightarrow X$ that minimizes the transportation cost
$\int_{x\in\RR^d} c(x, T(x)) d\mu(x)$ where $c$ is a cost function, under the constraint
that the measure over $X$ is the pushforward of $\mu$ by the map $T$, denoted by
$T_{\#}\mu$, which means
that for any measurable set $B$, $\int_{x\in T^{-1}(B)}d\mu(x) = \int_{y \in B} d\nu(y)$.
In the semi-discrete case, this means that for every site $x_i$, its pre-image
$T^{-1}(x_i)$ has the prescribed mass $p_i$: $\forall i, \int_{x \in T^{-1}(x_i)}d\mu(x) = p_i$.
In the remainder of this article, we will use the squared Euclidean distance 
for the cost: $c(x, y) = \norm{x - y}^2_2$.

As shown by Aurenhammer, Hoffmann and Aronov~\cite{aha}, any semi-discrete optimal 
transport map between a density and a discrete measure represented by a set of sites at 
location  $(x_i)_i$ with masses $(p_i)_i$ is such that the pre-images of the
$x_i$s through the optimal transport map correspond to the cells of a power diagram
(or more generally, a Laguerre diagram for arbitrary cost functions $c$) with a well-chosen 
set of weights.
For a power diagram with sites $(x_i)_i$ and associated weights $(\phi_i)_i$, we recall the expression of the power cell of site $x_i$:
$$Pow^\phi(x_i) = \set{x \in \RR^d \midd \forall 1 \le j \le N, \norm{x - x_i}^2 - \phi_i \le \norm{x - x_j}^2 - \phi_j}$$
Intuitively, a power diagram is a Voronoi diagram in which $x_i$'s cell size depends on the weight $\phi_i$. In particular, a higher weight value for $x_i$ relative to the weights of its neighbors $(x_j)$ leads to a larger cell. There is however no straightforward relationship between the value of a particular $\phi_i$ and its corresponding cell volume.

The condition $\sum_{i=1}^N p_i \delta_{x_i} = T_{\#}\mu$ means that the cells' measures need to respect the prescribed measure
on the discrete samples. Since there is no direct relationship between $(\phi_i)_i$ and volumes (or measures), we need to adjust the weights  $(\phi_i)_i$ through an energy minimization. 
Thus, the set of weights $(\phi_i)_i$ realizing the power diagram that corresponds to an optimal transport 
map can be expressed as the minimum of the following functional~\cite{aha}:

$$\Phi((\phi_i)_i) = \Sum_{i=1}^N \int_{x\in Pow^\phi(x_i)}(\norm{x - x_i} - \phi_i) d\mu(x) + \Sum_{i=1}^N p_i \phi_i$$

We will denote the transport map associated to the set of weights $\phi$ as $T_\phi$.

\section{Symmetrized semi-discrete optimal transport}\label{tosd}

\subsection{Goal and rationale} \label{goals}
 Given two measures $\mu$ and $\nu$ over $\mathbb{R}^d$ with compact supports, and two sets of samples with associated weights $((x_i,\phi_i))_i$ and $((y_i,\psi_i))_i$, 
 we denote $Pow^\phi$ (resp. $Pow^\psi$) the power diagram of $(x_i)_i$ (resp. $(y_i)_i$) associated with weights $(\phi_i)_i$ (resp. $(\psi_i)_i$) restricted to the support of $\nu$ (resp. $\mu$). In that context, we aim at finding the samples positions and weights such that $x_i$ is the barycenter of $Pow^\psi(y_i)$ and $y_i$ is the barycenter of $Pow^\phi(x_i)$ for all $i$, and power diagrams $Pow^\phi(x_i)$ and $Pow^\psi(y_i)$ respectively describing optimal transport maps from $\mu$ to $(x_i)_i$ and from $\nu$ to $(y,_i)_i$, both sets of samples being equipped with the discrete uniform measure. We denote these transport maps $T_\phi$ and $T_\psi$.

With such a construction, samples $(x_i)_i$ and $(y_i)_i$ are in one-to-one correspondence, and the power cells of the two semi-discrete optimal transport maps $Pow^\phi(x_i)$ and $Pow^\psi(y_i)$ are in correspondence as well. This will be exploited to design a mesh interpolation algorithm as shown in section~\ref{sec:displacement}.

This problem can be formalized as a constrained optimization problem. The
objective function is a combination of transport functionals that is minimal when power 
diagrams correspond to optimal transport maps. 
Constraints encode the fact that the sites $(x_i)$ and $(y_i)$ are located at
the barycenters of the respective cells $\Lagx$ and $\Lagy$.

We formulate it as follows:

\begin{align}
	\begin{split}
\min\limits_{(x_i, \phi_i, y_i, \psi_i)} & \Sum_{i=1}^N \int_{x\in Pow^\phi(x_i)}(\norm{x - x_i} - \phi_i) d\mu(x) + \Sum_{i=1}^N p_i \phi_i \\
+ & \Sum_{i=1}^N \int_{y\in Pow^\psi(y_i)}(\norm{y - y_i} - \psi_i) d\nu(x) + \Sum_{i=1}^N q_i \psi_i\\
\textrm{s.t. } & \int_{x\in Pow^\phi(x_i)} (x - y_i) d\mu(x) = 0\\
\textrm{and }&\int_{y\in Pow^\psi(y_i)} (y - x_i) d\nu(y) = 0
	\label{eq:objfun}
	\end{split}
\end{align}

By analogy with classical Centroidal Voronoi Tessellations construction algorithms such as Lloyd's algorithm,
we propose a fixed-point iteration to solve this constrained problem,
by alternatively minimizing the objective function and enforcing the constraints.

\subsection{Algorithm} \label{sec:mainalgo}
As input, our algorithm takes two measures $\mu$ and $\nu$ whose supports $Sp(\mu)$ and $Sp(\nu)$ are domains meshed with triangles (in 2-d) or tetrahedra (in 3-d). Measures $\mu$ and $\nu$ are defined as piecewise linear functions, and are entirely given by their values on the mesh vertices and linearly interpolated over triangles or tetrahedra.

Our algorithm starts by uniformly sampling the supports of both measures $\mu$ and $\nu$ following the method described by Levy and Bonneel~\cite{LB12} that samples each simplex proportionally to its area.
We then repeat until convergence the following operations. 

First, we optimize weights of the power diagram restricted to the support of $\mu$, $Pow^\phi(x_i)$, describing the semi-discrete transport map between $\mu$ and $(x_i)_i$ using standard semi-discrete optimal transport techniques~\cite{bl-semi-discrete}. 

We then move each sample $y_i$ to the barycenter $\bar{y}_i$ of the newly
computed power cell $Pow^\phi(x_i)$, accounting for measure $\mu$:  $\bar{y}_i = \frac{\int_{x\in\Lagx} x d\mu(x)}{\int_{x\in\Lagx}d\mu(x)}$. 

Then, we repeat the same operation by inverting the roles of the samples and measures -- computing the transport map between $\nu$ and $(y_i)_i$, and centering samples $(x_i)_i$ at the barycenter of $Pow^\psi(y_i)$ with respect to measure $\nu$.

\begin{algorithm}[!tbh]
  \begin{algorithmic}[1]
    \STATE $(x_i) := $ random sampling of $Sp(\nu)$ 
    \STATE $(y_i) := $ random sampling of $Sp(\mu)$
    \WHILE{not converged}
    \STATE $(\phi_i) :=$ semi-discrete optimal transport weights from $\mu$ to $(x_i)$
    \STATE $(y_i) :=$ centroids of $Pow^\phi(x_i)$
    \STATE $(\psi_i) :=$ semi-discrete optimal transport weights from $\nu$ to $(y_i)$
    \STATE $(x_i) :=$ centroids of $Pow^\psi(y_i)$
    \ENDWHILE
  \end{algorithmic}
  \caption{Alternating algorithm for computing symetrized optimal transport maps}
  \label{alg:seq}
\end{algorithm}

\begin{figure}[!tbh]
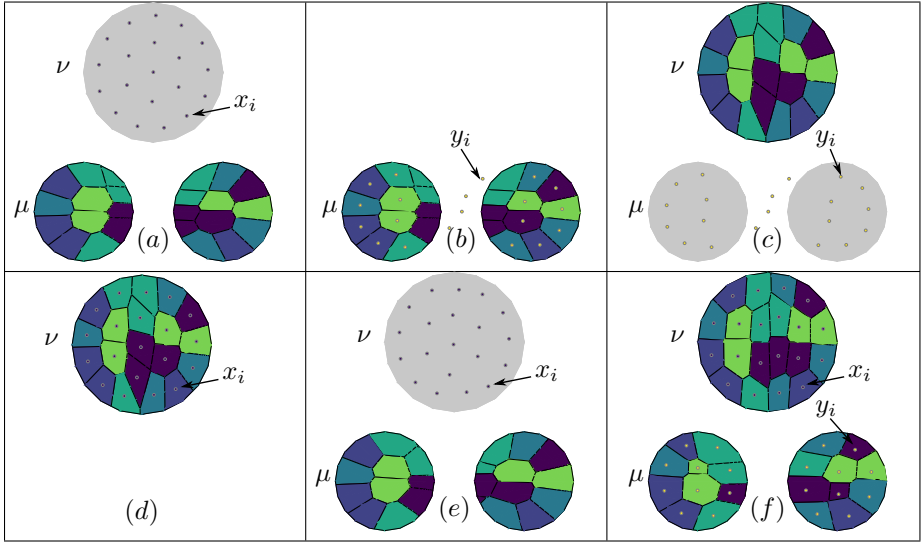

  \centering
  \begin{tabular}{|c|c|c|}
    \hline 
    \includesvg[scale=0.35]{algo-alterne/step0.svg}&
    \includesvg[scale=0.35]{algo-alterne/step1.svg}&
    \includesvg[scale=0.35]{algo-alterne/step2.svg}\\
    \hline
    \includesvg[scale=0.35]{algo-alterne/step3.svg}&
    \includesvg[scale=0.35]{algo-alterne/step4.svg}&
    \includesvg[scale=0.35]{algo-alterne/final.svg}\\
    \hline
    \end{tabular}
  \caption{First steps of the algorithm and final setting. Steps $(a-d)$ correspond to a single iteration of Alg.~\protect\ref{alg:seq}. 
  $(a)$ Measure $\mu$ is transported to samples $(x_i)$, located on the support of $\nu$.
  $(b)$ Samples $(y_i)$ are relocated to the barycenters of the new power cells.
  $(c)$ Measure $\nu$ is transported to samples $(y_i)$, located on the support of $\mu$.
  $(d)$ Samples $(x_i)$ are relocated to the barycenters of the new power cells.
  $(e)$ Transport from $\mu$ to $(x_i)$ for the second iteration. 
  $(f)$ Transport maps and samples at convergence after several iterations.}
  \label{example-algo}
\end{figure}

We show that even this seemingly simple algorithm that symmetrizes the notion of semi-discrete optimal transport produces displacement interpolation results that well capture discontinuous behavior in the transport maps.
A sample run of Algorithm~\ref{alg:seq} on a 2-d example is shown in Fig.~\ref{example-algo}.

Each semi-discrete optimal transport computation results in an optimization, typically performed using an iterative solver (L-BFGS in our case). 
For the two optimal transport optimizations performed at the first (outer) iteration, in practice we initialize transport weights with a constant value, resulting in Voronoi diagrams. For the remaining iterations we employ a warm restart strategy: the optimized values of $(\phi_i)_i$ and $(\psi_i)_i$ from the previous iteration are reused as initial guesses. We repeat these iterations a fixed number of times.
We found that 100 iterations were enough in practice to reach convergence in all our interpolation examples.

\section{Mesh interpolation algorithm} \label{sec:displacement}

The main application to our symmetrized semi-discrete optimal transport algorithm 
is displacement interpolation -- or warping -- between shapes, typically in 2 or 3 dimensions.

In a nutshell, our alternating algorithm leads to cells $Pow^\phi(x_i)$ and
$Pow^\psi(y_i)$ having the same geometry for all $i$ (with the possible exception
of cells on the mesh boundary).
This allows us to come up with a simple interpolation technique: since
corresponding power cells have the same geometry, we can put in
their vertices correspondence and linearly interpolate between them.
In addition, power cells boundaries tend to be aligned with the transport maps discontinuities, which helps capturing tearing during displacement interpolations. 
This allows to define the following 2-d mesh interpolation algorithm.

Our algorithm is interested in restricted power diagrams.
We will use the classification used by Nivoliers~\cite{nivoliers-vsdm} to characterize the vertices of such meshes. A vertex of a restricted power diagram is necessarily of one of the three following types:

\begin{itemize}
\item \emph{type i} a vertex that originates from the underlying mesh and that does not depend on
  the power diagram,
\item \emph{type ii} a vertex that is located at the intersection of an edge separating 
  two power cells, and an edge of the underlying mesh,
\end{itemize}

\begin{itemize}
\item \emph{type iii} a vertex that is at the intersection between three cells of the power diagram.
  Such a vertex can be uniquely identified by the triplet $(i,j,k)$ of surrounding
  cells indices.
\end{itemize}

\begin{wrapfigure}{r}{0.5\textwidth}
  \centering
  \includesvg[width=0.4\textwidth]{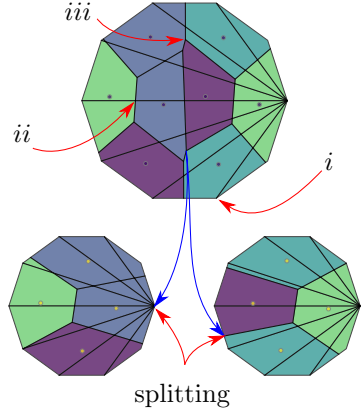}
  \caption{Vertices with type $i$, $ii$ and $iii$, and a vertex that splits during the interpolation.}
  \label{fig:vertices-types}
\end{wrapfigure}

As illustrated in Figure~\ref{fig:vertices-types}, we observe that,
in practice, corresponding power cells have similar shapes.
More formally, for given indices $i$ and $j$, we found that whenever the cells $Pow^\phi(x_i)$ and 
$Pow^\phi(x_j)$ are neighbors, $Pow^\psi(y_i)$ and $Pow^\psi(y_j)$ are also
neighbors.
While we do not provide a formal proof for this, we observed this behavior in almost all of our examples.
In the rare situations where this is not the case, adding more iterations to our algorithm resolves these instances, which allows us to handle these cases.

Each power diagram vertex being uniquely defined by the power cells it is incident to, and power cells being in one-to-one correspondence, it is easy to check whether a vertex of one power diagram can be matched with a vertex of the other, and match them. 

When intersecting the power diagram with a mesh, this translates as
the fact that there is at most a one-to-one mapping of type iii vertices,
whenever corresponding vertices exist in both restricted power diagrams.
We thus identify cell vertices based on their neighboring cells, and
linearly interpolate between the two vertices that share the same neighboring
cells in the two domains.

Correspondence between type ii vertices that lie on the boundary follows the
same principle. If such a vertex is adjacent to cells $Pow^\phi(x_i)$ and $Pow^\phi(x_j)$, it is identified by the triplet $(i, j, 0)$, where $0$ represents the
fact that the ``third cell'' it is adjacent to is in fact the outside of the mesh. 

However, whenever the meshes present topological discontinuities, some type iii
vertices are bound to split during the transport, typically giving birth to two
type ii vertices. In these cases, we need to duplicate the type iii vertices identified
by $(i, j, k)$ into two new type ii vertices among $(i, j, 0)$, $(i, 0, k)$ and $(0, j, k)$,
and seek corresponding vertices on the other side.
A similar treatement is applied to type ii vertices that split into type i vertices.

The fact that there is no one-to-one correspondence between vertices involved in
discontinuities offers us a practical criterion for identifying cells vertices
that lie on a discontinuity.

At last, we need to account for type i vertices lying on the boundary.
In the best of cases, there is only one vertex identified by $(i, 0, 0)$ in each
of the restricted Voronoi diagrams, and we can associate them right away.
In the worst case, there are several vertices represented by the same identifier
on each side. In this situation, we chose to associate all the equivalent vertices
on one restricted Voronoi diagram to a single vertex of the opposite. This potentially
causes some artefacts in the morph.

\section{Results} \label{results}

Computation of the transport maps is done using the \texttt{HLBFGS} library~\cite{hlbfgs}, and the algorithm has been implemented within the Graphite library~\cite{graphite}. 

\subsection{Symmetrized transport computation}

The overall complexity of our algorithm is dominated by the optimal transport computation, since the recentering complexity is negligible.
Requiring multiple calls to optimal transport optimizations makes the overall procedure 
relatively costly -- though of performance similar to iterative semi-discrete
optimal transport computations of fluid dynamics~\cite{gallouet2018lagrangian}
that perform similar iterations. Due to our warm restart, transport map computations are typically much faster after the first iteration. In practice, in all our 2-d and 3-d examples, the entire process takes approximately 7--8 minutes for 200 samples and 27-37 minutes for 10k samples, on an Intel Xeon E5-1650 6-core machine at 3.5GHz.

\subsection{Interpolations in 2-d}\label{sec:2d}

\begin{figure*}[!tbh]
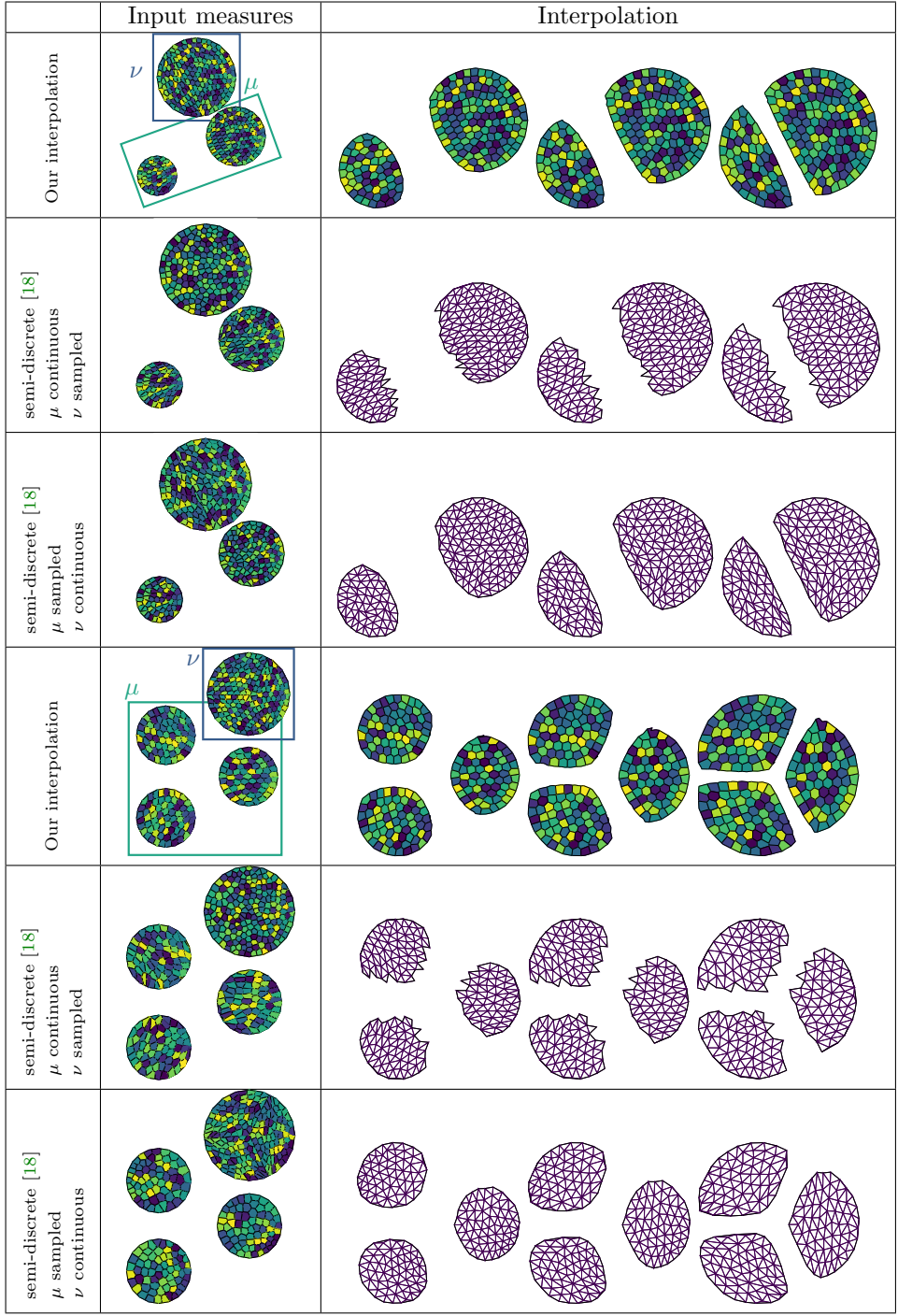

    \centering
    \begin{tabular}{|c|c|c|c|c|}
    \hline
    &Input measures & Interpolation  \\       
    \hline   
      \rotatebox[origin=l]{90}{\scriptsize{Our interpolation} } &
    \includesvg[width=2.7cm]{interpolation/2d/shape2/shape2-scene.svg} &
    \includesvg[height=2.0cm]{interpolation/2d/shape2/shape2-step0.svg}
    \includesvg[height=2.0cm]{interpolation/2d/shape2/shape2-step1.svg}
    \includesvg[height=2.0cm]{interpolation/2d/shape2/shape2-step2.svg} \\	
    \hline
     \rotatebox[origin=l]{90}{\scriptsize{semi-discrete~\cite{bl-semi-discrete}}}
      \rotatebox[origin=l]{90}{\scriptsize{$\mu$ continuous}} \rotatebox[origin=l]{90}{\scriptsize{$\nu$ sampled}} &
    \includesvg[width=2.7cm]{interpolation/2d/shape2/shape2-nonsym.svg} & 
    \includesvg[height=2.0cm]{interpolation/2d/shape2/shape2-ns-step2.svg}
    \includesvg[height=2.0cm]{interpolation/2d/shape2/shape2-ns-step1.svg}
    \includesvg[height=2.0cm]{interpolation/2d/shape2/shape2-ns-step0.svg}\\

    \hline
     \rotatebox[origin=l]{90}{\scriptsize{semi-discrete~\cite{bl-semi-discrete}}}
      \rotatebox[origin=l]{90}{\scriptsize{$\mu$ sampled}} \rotatebox[origin=l]{90}{\scriptsize{$\nu$ continuous}} &
    \includesvg[width=2.7cm]{interpolation/2d/shape2/shape2-nonsym-rev.svg}&
    \includesvg[height=2.0cm]{interpolation/2d/shape2/shape2-ns-rev-step0.svg}
    \includesvg[height=2.0cm]{interpolation/2d/shape2/shape2-ns-rev-step1.svg}
    \includesvg[height=2.0cm]{interpolation/2d/shape2/shape2-ns-rev-step2.svg}\\
	
	\hline
	\rotatebox[origin=l]{90}{\scriptsize{Our interpolation} } &
	\includesvg[width=2.4cm]{interpolation/2d/shape4/shape4-scene.svg} &
    \includesvg[height=2.3cm]{interpolation/2d/shape4/shape4-step0.svg}
    \includesvg[height=2.3cm]{interpolation/2d/shape4/shape4-step1.svg}
    \includesvg[height=2.3cm]{interpolation/2d/shape4/shape4-step2.svg} \\ 
	
	\hline
     \rotatebox[origin=l]{90}{\scriptsize{semi-discrete~\cite{bl-semi-discrete}}}
      \rotatebox[origin=l]{90}{\scriptsize{$\mu$ continuous}} \rotatebox[origin=l]{90}{\scriptsize{$\nu$ sampled}} &
    \includesvg[width=2.4cm]{interpolation/2d/shape4/shape4-nonsym.svg} & 
    \includesvg[height=2.3cm]{interpolation/2d/shape4/shape4-ns-step2.svg}
    \includesvg[height=2.3cm]{interpolation/2d/shape4/shape4-ns-step1.svg}
    \includesvg[height=2.3cm]{interpolation/2d/shape4/shape4-ns-step0.svg}\\	
	
	\hline
     \rotatebox[origin=l]{90}{\scriptsize{semi-discrete~\cite{bl-semi-discrete}}}
      \rotatebox[origin=l]{90}{\scriptsize{$\mu$ sampled}} \rotatebox[origin=l]{90}{\scriptsize{$\nu$ continuous}} &	
	    \includesvg[width=2.4cm]{interpolation/2d/shape4/shape4-nonsym-rev.svg}&		
    \includesvg[height=2.3cm]{interpolation/2d/shape4/shape4-ns-rev-step0.svg}
    \includesvg[height=2.3cm]{interpolation/2d/shape4/shape4-ns-rev-step1.svg}
    \includesvg[height=2.3cm]{interpolation/2d/shape4/shape4-ns-rev-step2.svg}\\	
    \hline
    \end{tabular}
    \caption{Interpolations between a disk and two disks, and between a disk and three disks, using 200 samples. 
      The non-symmetric approach poorly approximates discontinuities when the single 
	  disk is sampled. Our method best captures discontinuities in all cases.}\label{fig:shape2}
\end{figure*}

\begin{figure}[!h]
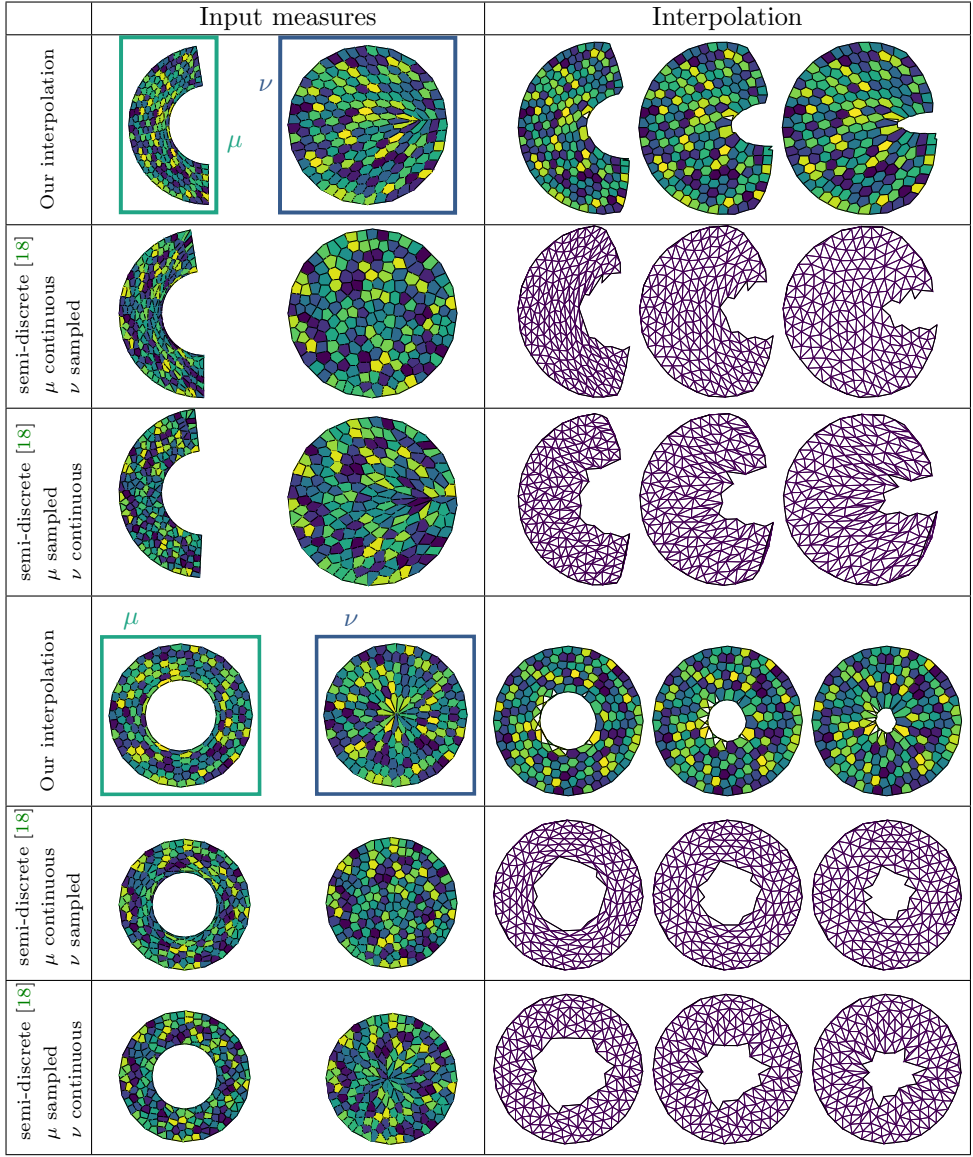

    \centering
    \begin{tabular}{|c|c|c|}
      \hline
      &Input measures & Interpolation \\
      \hline   
      \rotatebox[origin=l]{90}{\scriptsize{Our interpolation} } &
      \includesvg[width=4.5cm]{interpolation/2d/shape6/shape6-scene.svg} &
      \includesvg[height=2.3cm]{interpolation/2d/shape6/shape6-step0.svg}
      \includesvg[height=2.3cm]{interpolation/2d/shape6/shape6-step1.svg}
      \includesvg[height=2.3cm]{interpolation/2d/shape6/shape6-step2.svg} \\ 
      \hline
     \rotatebox[origin=l]{90}{\scriptsize{semi-discrete~\cite{bl-semi-discrete}}}
     \rotatebox[origin=l]{90}{\scriptsize{$\mu$ continuous}}
     \rotatebox[origin=l]{90}{\scriptsize{$\nu$ sampled}} &
      \includesvg[width=4.5cm]{interpolation/2d/shape6/shape6-nonsym.svg} & 
      \includesvg[height=2.3cm]{interpolation/2d/shape6/shape6-ns-step2.svg}
      \includesvg[height=2.3cm]{interpolation/2d/shape6/shape6-ns-step1.svg}
      \includesvg[height=2.3cm]{interpolation/2d/shape6/shape6-ns-step0.svg}\\
      \hline
     \rotatebox[origin=l]{90}{\scriptsize{semi-discrete~\cite{bl-semi-discrete}}}
     \rotatebox[origin=l]{90}{\scriptsize{$\mu$ sampled}}
     \rotatebox[origin=l]{90}{\scriptsize{$\nu$ continuous}} &	
      \includesvg[width=4.5cm]{interpolation/2d/shape6/shape6-nonsym-rev.svg}&
      \includesvg[height=2.3cm]{interpolation/2d/shape6/shape6-ns-rev-step0.svg}
      \includesvg[height=2.3cm]{interpolation/2d/shape6/shape6-ns-rev-step1.svg}
      \includesvg[height=2.3cm]{interpolation/2d/shape6/shape6-ns-rev-step2.svg}\\
      \hline
      \rotatebox[origin=l]{90}{\scriptsize{Our interpolation} } &
    \includesvg[width=5cm]{interpolation/2d/shape12/shape12-scene.svg} &
    \includesvg[width=2cm]{interpolation/2d/shape12/shape12-step0.svg}
    \includesvg[width=2cm]{interpolation/2d/shape12/shape12-step1.svg}
    \includesvg[width=2cm]{interpolation/2d/shape12/shape12-step2.svg} \\ 
      \hline
     \rotatebox[origin=l]{90}{\scriptsize{semi-discrete~\cite{bl-semi-discrete}}}
     \rotatebox[origin=l]{90}{\scriptsize{$\mu$ continuous}}
     \rotatebox[origin=l]{90}{\scriptsize{$\nu$ sampled}} &
    \includesvg[width=4.5cm]{interpolation/2d/shape12/shape12-nonsym.svg} & 
    \includesvg[width=2cm]{interpolation/2d/shape12/shape12-ns-step2.svg}
    \includesvg[width=2cm]{interpolation/2d/shape12/shape12-ns-step1.svg}
    \includesvg[width=2cm]{interpolation/2d/shape12/shape12-ns-step0.svg}\\
      \hline
     \rotatebox[origin=l]{90}{\scriptsize{semi-discrete~\cite{bl-semi-discrete}}}
     \rotatebox[origin=l]{90}{\scriptsize{$\mu$ sampled}}
     \rotatebox[origin=l]{90}{\scriptsize{$\nu$ continuous}} &	
    \includesvg[width=4.5cm]{interpolation/2d/shape12/shape12-nonsym-rev.svg}&
    \includesvg[width=2cm]{interpolation/2d/shape12/shape12-ns-rev-step0.svg}
    \includesvg[width=2cm]{interpolation/2d/shape12/shape12-ns-rev-step1.svg}
    \includesvg[width=2cm]{interpolation/2d/shape12/shape12-ns-rev-step2.svg}\\
      \hline	  
    \end{tabular}
    \caption{Interpolation between a disk and non-convex shapes. A
      discontinuity appears on the disk to create the non-convexity.}\label{fig:shape6}
\end{figure}

\begin{figure}[!tbh]
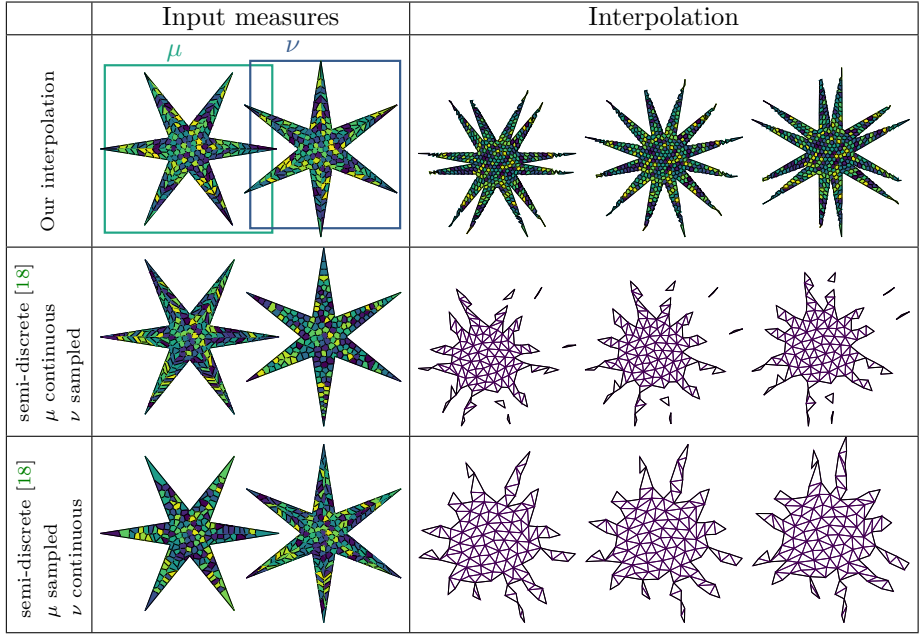

    \centering
    \begin{tabular}{|c|c|c|}
    \hline
      & Input measures & Interpolation \\
    \hline   
      \rotatebox[origin=l]{90}{\scriptsize{Our interpolation} } &
    \includesvg[width=4cm]{interpolation/2d/shape7/shape7-scene.svg} &
    \includesvg[width=2.1cm]{interpolation/2d/shape7/shape7-step0.svg}
    \includesvg[width=2.1cm]{interpolation/2d/shape7/shape7-step1.svg}
    \includesvg[width=2.1cm]{interpolation/2d/shape7/shape7-step2.svg} \\ 
    \hline
     \rotatebox[origin=l]{90}{\scriptsize{semi-discrete~\cite{bl-semi-discrete}}}
     \rotatebox[origin=l]{90}{\scriptsize{$\mu$ continuous}}
     \rotatebox[origin=l]{90}{\scriptsize{$\nu$ sampled}} &
    \includesvg[width=4cm]{interpolation/2d/shape7/shape7-nonsym.svg} & 
    \includesvg[width=2.1cm]{interpolation/2d/shape7/shape7-ns-step2.svg}
    \includesvg[width=2.1cm]{interpolation/2d/shape7/shape7-ns-step1.svg}
    \includesvg[width=2.1cm]{interpolation/2d/shape7/shape7-ns-step0.svg}\\
    \hline
     \rotatebox[origin=l]{90}{\scriptsize{semi-discrete~\cite{bl-semi-discrete}}}
     \rotatebox[origin=l]{90}{\scriptsize{$\mu$ sampled}}
     \rotatebox[origin=l]{90}{\scriptsize{$\nu$ continuous}} &
    \includesvg[width=4cm]{interpolation/2d/shape7/shape7-nonsym-rev.svg}&
    \includesvg[width=2.1cm]{interpolation/2d/shape7/shape7-ns-rev-step0.svg}
    \includesvg[width=2.1cm]{interpolation/2d/shape7/shape7-ns-rev-step1.svg}
    \includesvg[width=2.1cm]{interpolation/2d/shape7/shape7-ns-rev-step2.svg}\\

    \hline
    \end{tabular}
    \caption{Interpolation between a star and another star with branches
      equidistant to the first star's branches. Each branch splits in two equal
      parts that join with their neighbors to create the target branches.}
    \label{fig:shape7}
\end{figure}

\begin{figure}[!h]
  \centering
  \begin{tabular}{|c|c|c|c|}
    \hline
    & Source & Interpolation & Target \\
    & distribution & & distribution \\
    \hline
      \rotatebox[origin=l]{90}{\scriptsize{Our interpolation} } &  \includegraphics[width=2.10cm]{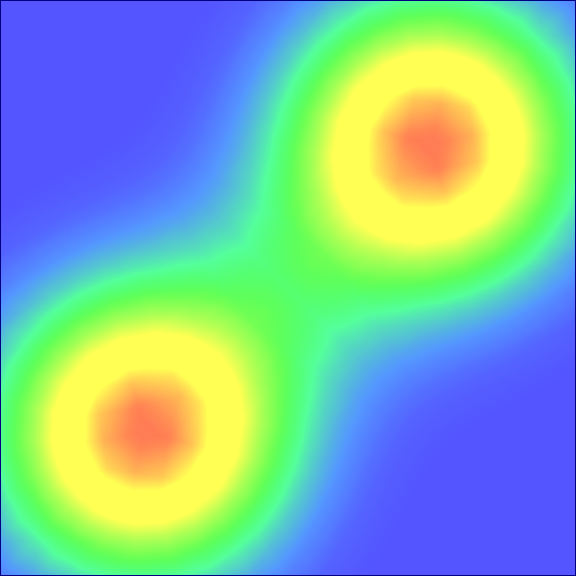} &
  \includesvg[width=2.10cm]{interpolation/2d/shape0/shape0-step05.svg}
  \includesvg[width=2.10cm]{interpolation/2d/shape0/shape0-step10.svg}
  \includesvg[width=2.10cm]{interpolation/2d/shape0/shape0-step15.svg} &
                                                                                 \includegraphics[width=2.10cm]{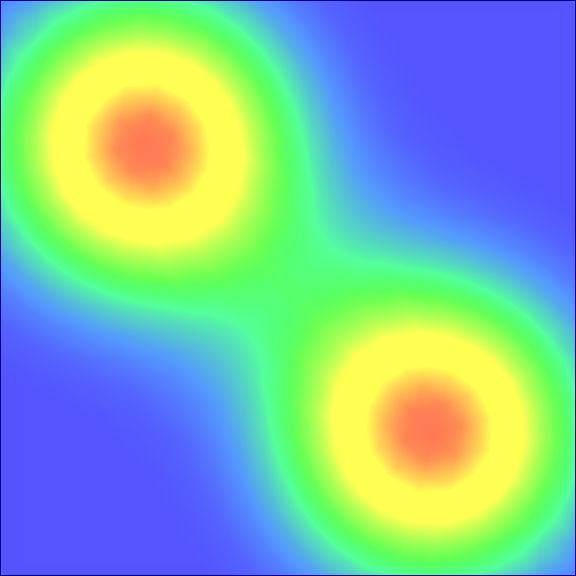} \\
    \hline
         \rotatebox[origin=l]{90}{\scriptsize{semi-discrete~\cite{bl-semi-discrete}}}
         \rotatebox[origin=l]{90}{\scriptsize{$\mu$ continuous}}
         \rotatebox[origin=l]{90}{\scriptsize{$\nu$ sampled}} &
    \includegraphics[width=2.10cm]{interpolation/2d/shape0/shape0-initial.png} &
  \includesvg[width=2.10cm]{interpolation/2d/shape0/shape0-ns-step0.svg}
  \includesvg[width=2.10cm]{interpolation/2d/shape0/shape0-ns-step1.svg}
  \includesvg[width=2.10cm]{interpolation/2d/shape0/shape0-ns-step2.svg} &

                                   \includegraphics[width=2.10cm]{interpolation/2d/shape0/shape0-final.png} \\
    \hline
     \rotatebox[origin=l]{90}{\scriptsize{semi-discrete~\cite{bl-semi-discrete}}}
     \rotatebox[origin=l]{90}{\scriptsize{$\mu$ sampled}}
     \rotatebox[origin=l]{90}{\scriptsize{$\nu$ continuous}} &
    \includegraphics[width=2.10cm]{interpolation/2d/shape0/shape0-initial.png} &
  \includesvg[width=2.10cm]{interpolation/2d/shape0/shape0-ns-rev-step2.svg} 
  \includesvg[width=2.10cm]{interpolation/2d/shape0/shape0-ns-rev-step1.svg}
  \includesvg[width=2.10cm]{interpolation/2d/shape0/shape0-ns-rev-step0.svg}&
                                                                                                                                                                        \includegraphics[width=2.10cm]{interpolation/2d/shape0/shape0-final.png} \\
    \hline
  \end{tabular}
  \caption{Interpolation between two squares equipped with non-uniform densities, each consisting of the sum of two gaussians}
  \label{fig:shape0}
\end{figure}

We compare our symmetrized algorithm with classical semi-discrete optimal transport~\cite{bl-semi-discrete} for 2D interpolation in Figures~\ref{fig:shape9}, \ref{fig:shape2}, \ref{fig:shape6}, \ref{fig:shape7} and \ref{fig:shape0}.

\paragraph{Measures with uniform densities}

Figures~\ref{fig:shape9} and ~\ref{fig:shape2} illustrate our algorithm on sets of disks -- two disks interpolated against two other disks, one disk interpolated against two disks, and one disk interpolated against three disks -- and compare it with the (non-symmetric) semi-discrete approach of Lévy~\cite{bl-semi-discrete}. 
When a single disk is interpolated with a shape consisting of two or three disks, the classical semi-discrete approach works well only when considering the single disk as the continuous measure and approximating the other (non-connected) shape with samples. However, appropriately choosing the continuous measure is not possible when interpolating between two non-connected shapes: in that case, our method still nicely captures tearing (Fig.~\ref{fig:shape9}) while a non-symmetric approach poorly approximates it.

Figure~\ref{fig:shape6} shows an interpolation from a single disk to connected but non-convex negatively-curved shapes, thus resulting in discontinuous transport maps by construction~\cite{chodosh2015discontinuity}. The non-symmetrized algorithm fails at capturing transport map discontinuities, both when transporting from or to the single disk, while our symmetric approach captures them well.

Figure~\ref{fig:shape7} shows a more complex interpolation between stars, where each branch separates into two equal parts during interpolation. Our approach better preserves the thin branch structures.

\setlength{\extrarowheight}{0pt}
\setlength{\tabcolsep}{3pt}

\paragraph{Measures with non-uniform densities}

As illustrated in Figure~\ref{fig:shape0}, our algorithm also handles measures that are not uniform over their support. We observe that
the measures are interpolated in a consistent way, with the two peaks from
the source distribution splitting and joining into the two peaks from the target
distribution. In contrast, the non-symmetrized algorithm results in interpolation of near uniform densities in this example.

\subsection{Interpolations in 3-d}
We demonstrate our algorithm on 3-d examples.

\begin{figure}[!h]
  \centering
  \begin{tabular}{|c|c|c|}
    \hline
     & Input measures & Interpolation \\
    \hline
    \rotatebox[origin=l]{90}{\scriptsize{Our interpolation} } &      \includegraphics[scale=0.08]{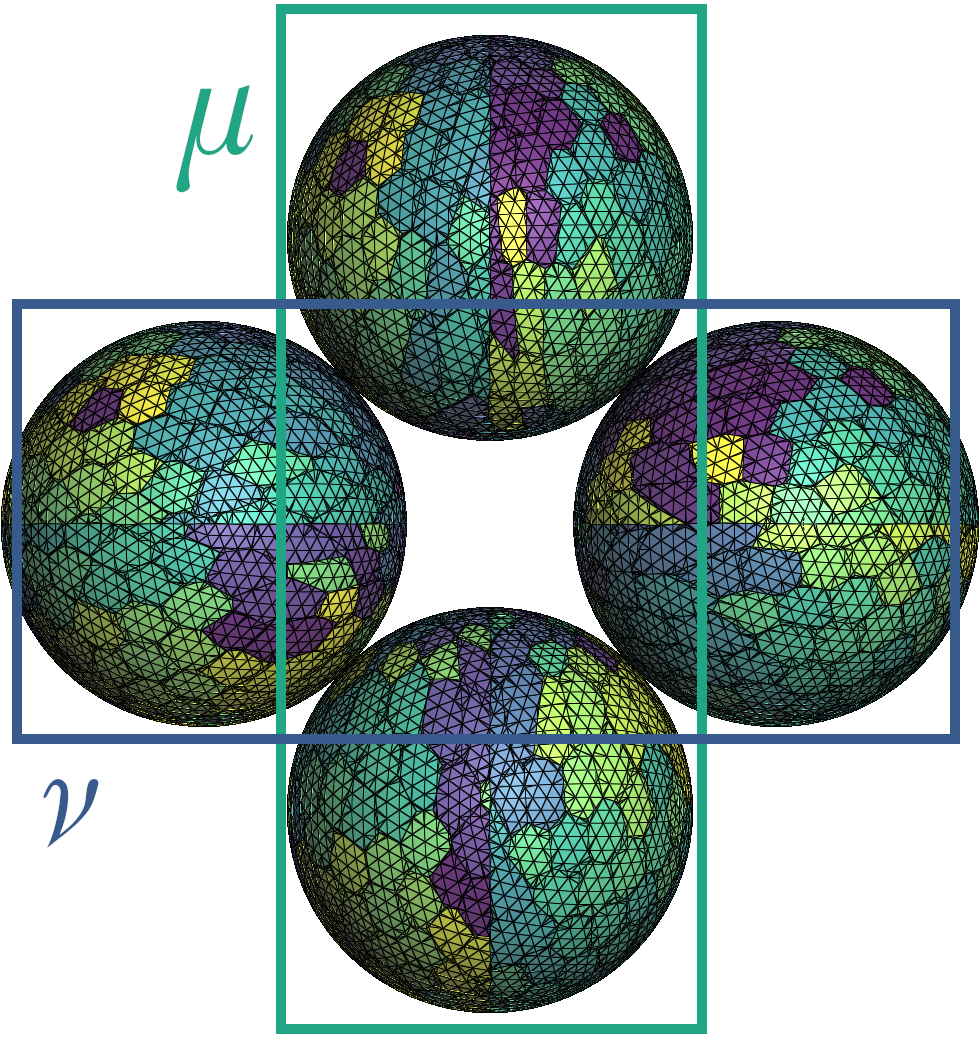} &
    \includegraphics[scale=0.08]{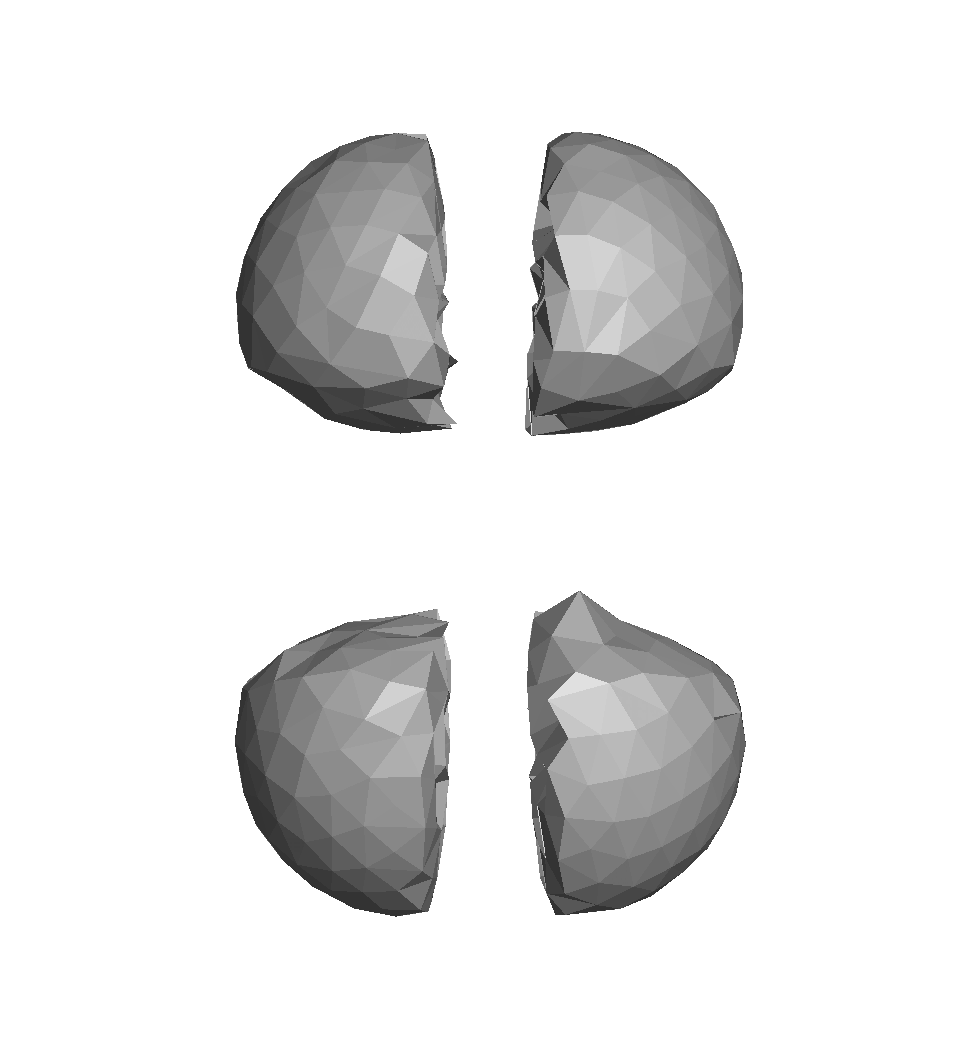}
    \includegraphics[scale=0.08]{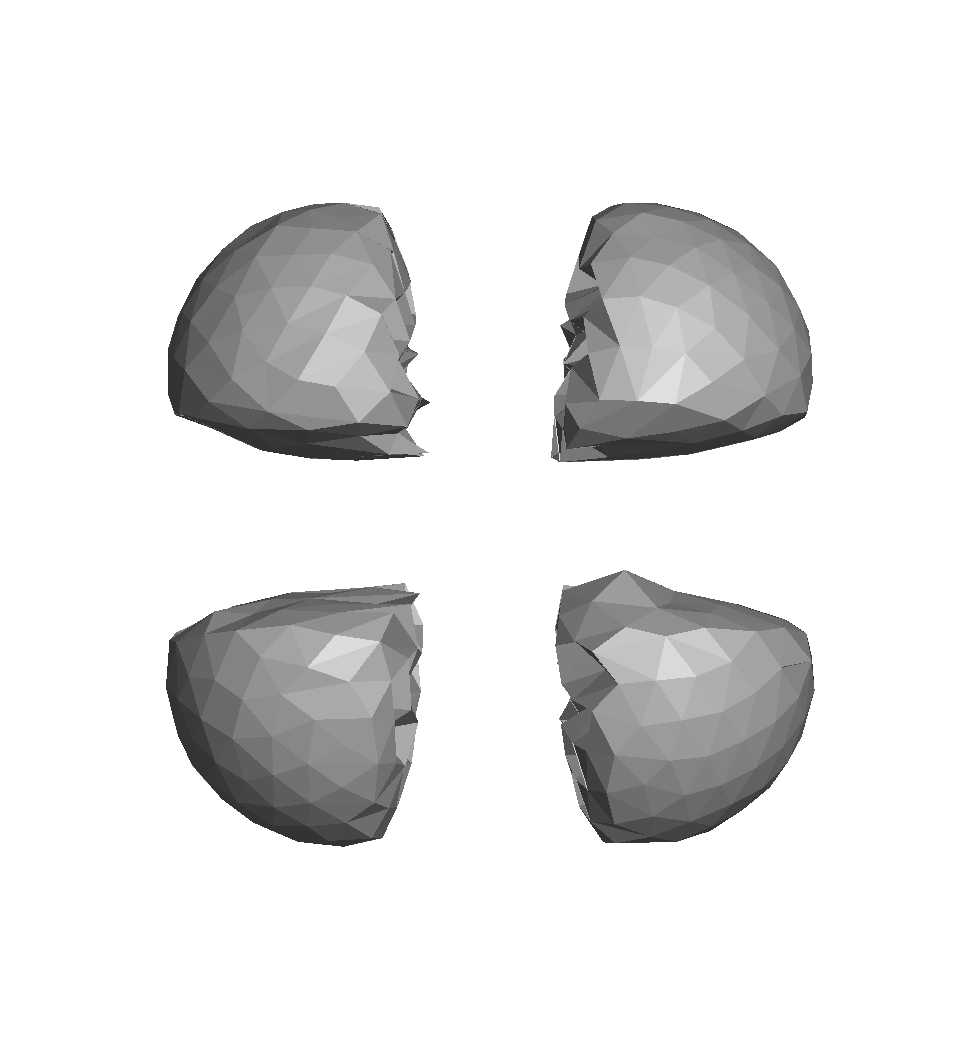}
    \includegraphics[scale=0.08]{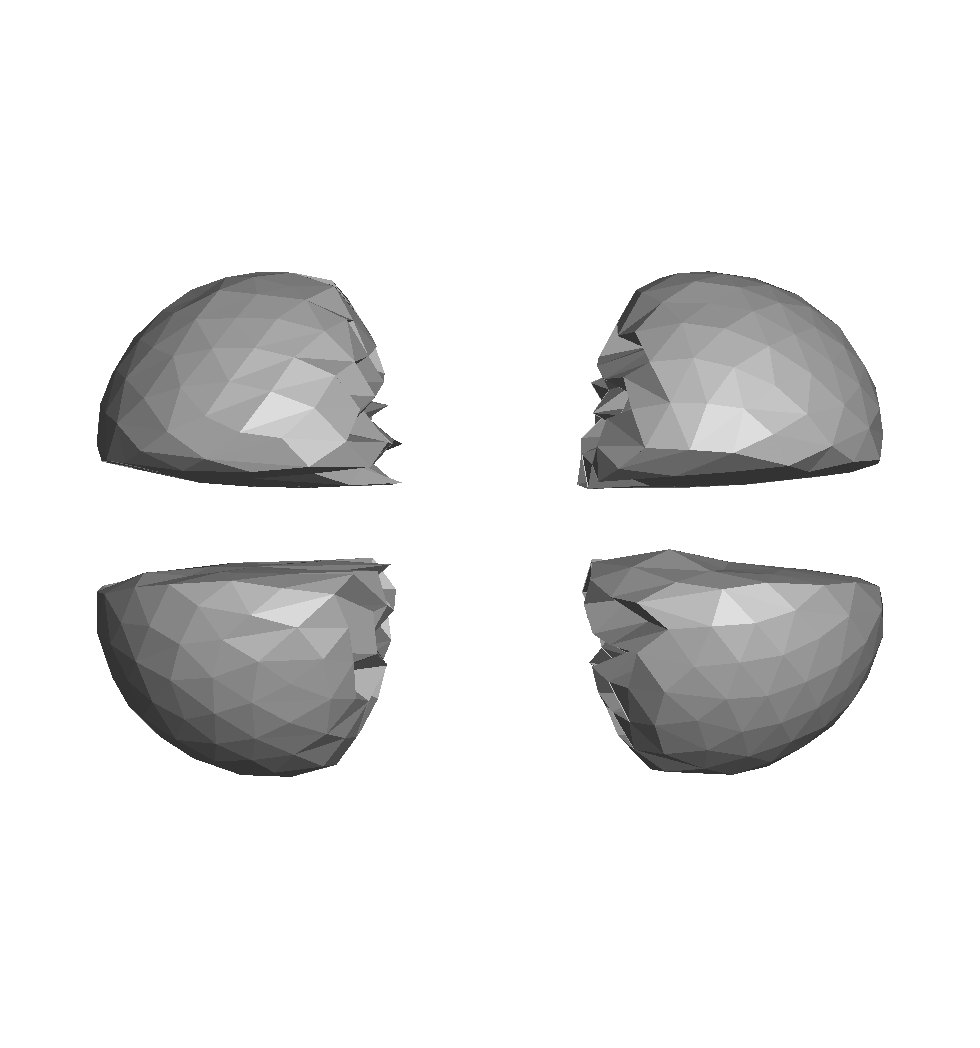} \\
    \hline
    \rotatebox[origin=l]{90}{\scriptsize{semi-discrete~\cite{bl-semi-discrete}}}
    \rotatebox[origin=l]{90}{\scriptsize{$\mu$ continuous}}
    \rotatebox[origin=l]{90}{\scriptsize{$\nu$ sampled}} &
    \includegraphics[scale=0.08]{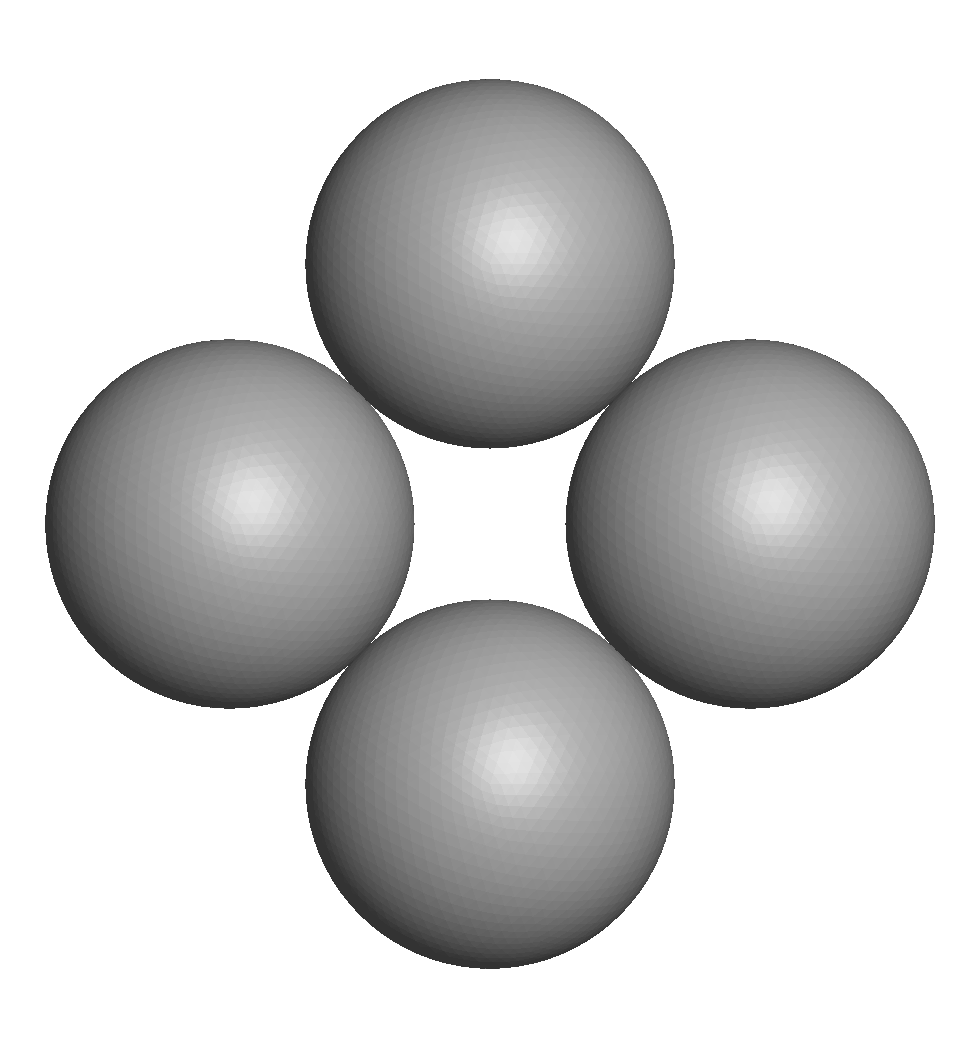} &
    \includegraphics[scale=0.08]{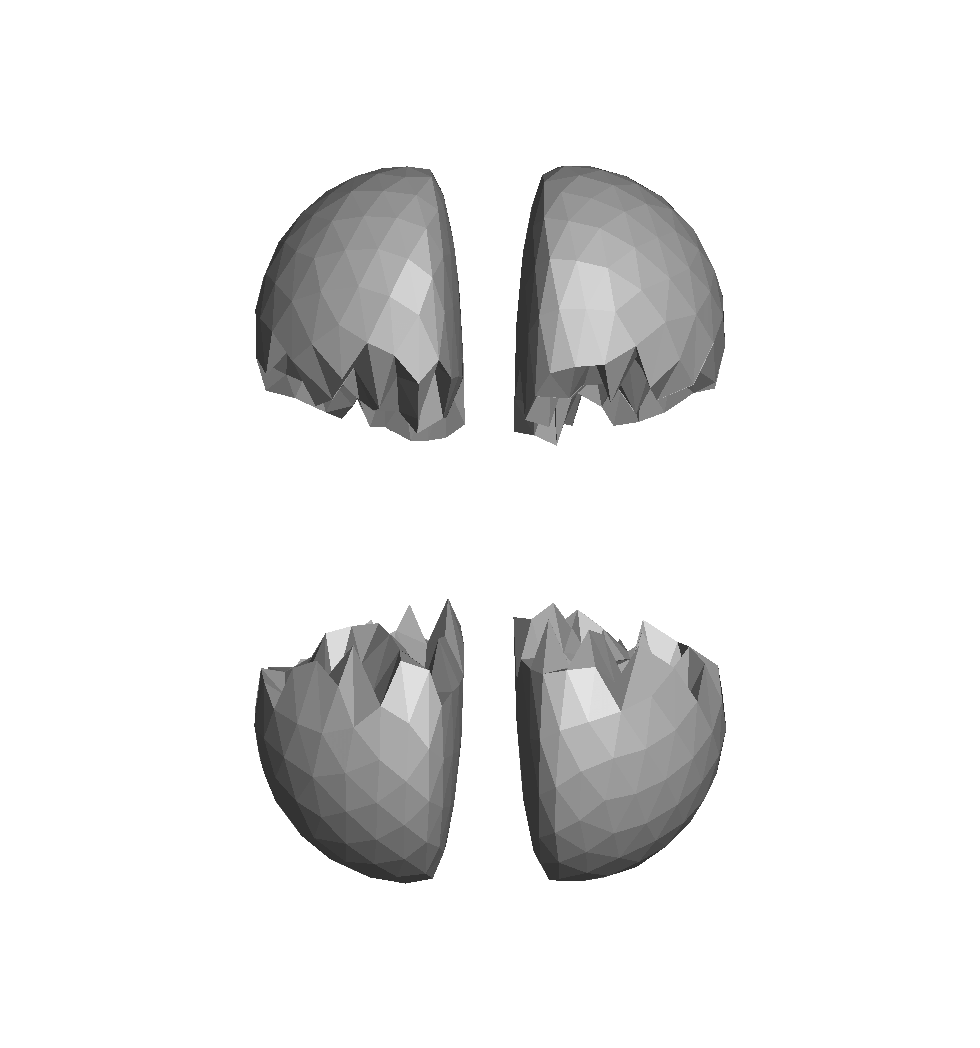}
    \includegraphics[scale=0.08]{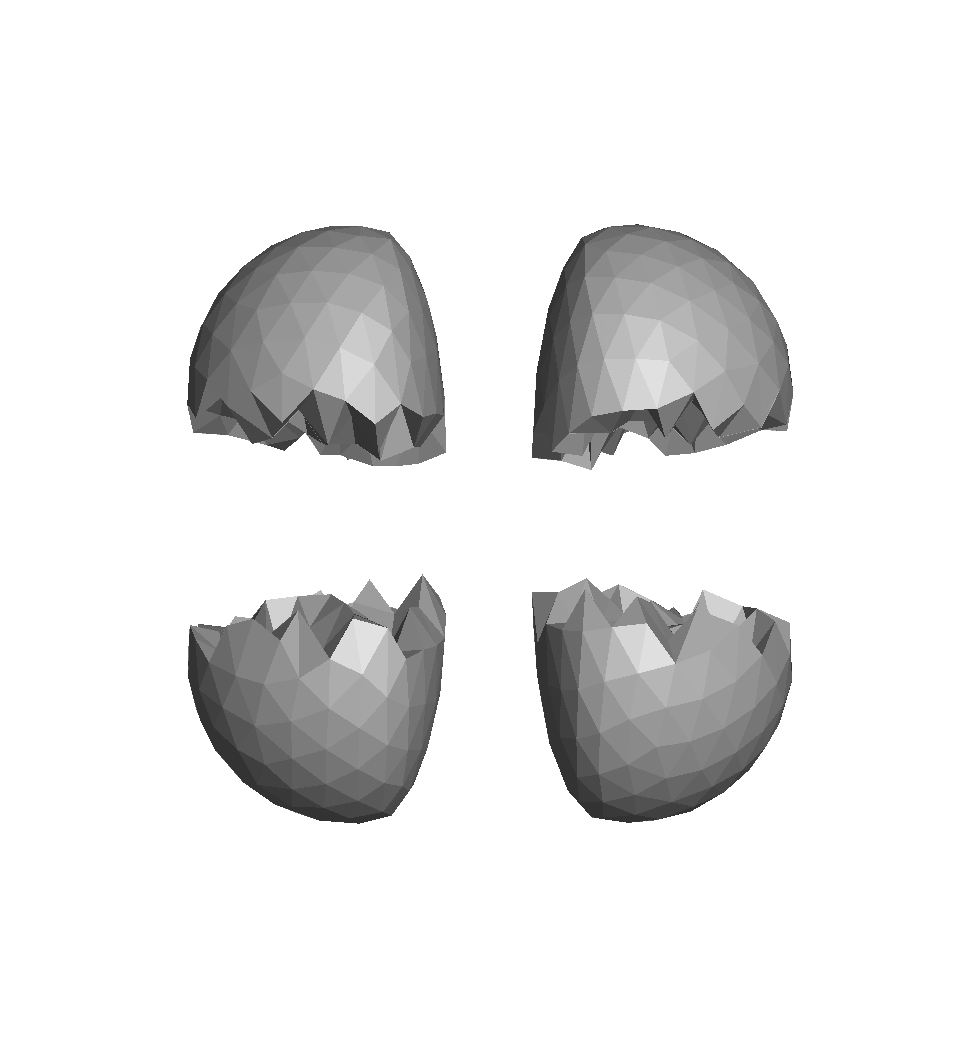}
    \includegraphics[scale=0.08]{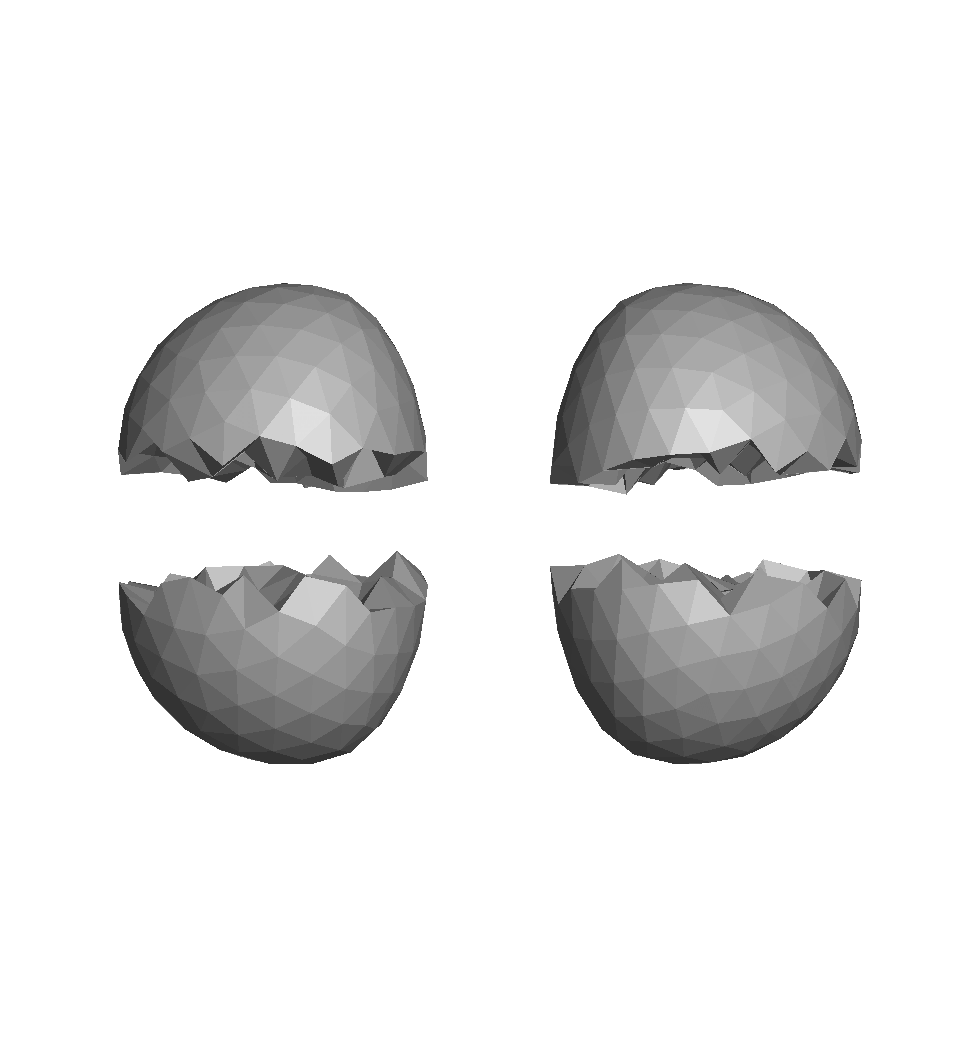} \\
    \hline
     \rotatebox[origin=l]{90}{\scriptsize{semi-discrete~\cite{bl-semi-discrete}}}
     \rotatebox[origin=l]{90}{\scriptsize{$\mu$ sampled}}
     \rotatebox[origin=l]{90}{\scriptsize{$\nu$ continuous}} &    
    \includegraphics[scale=0.08]{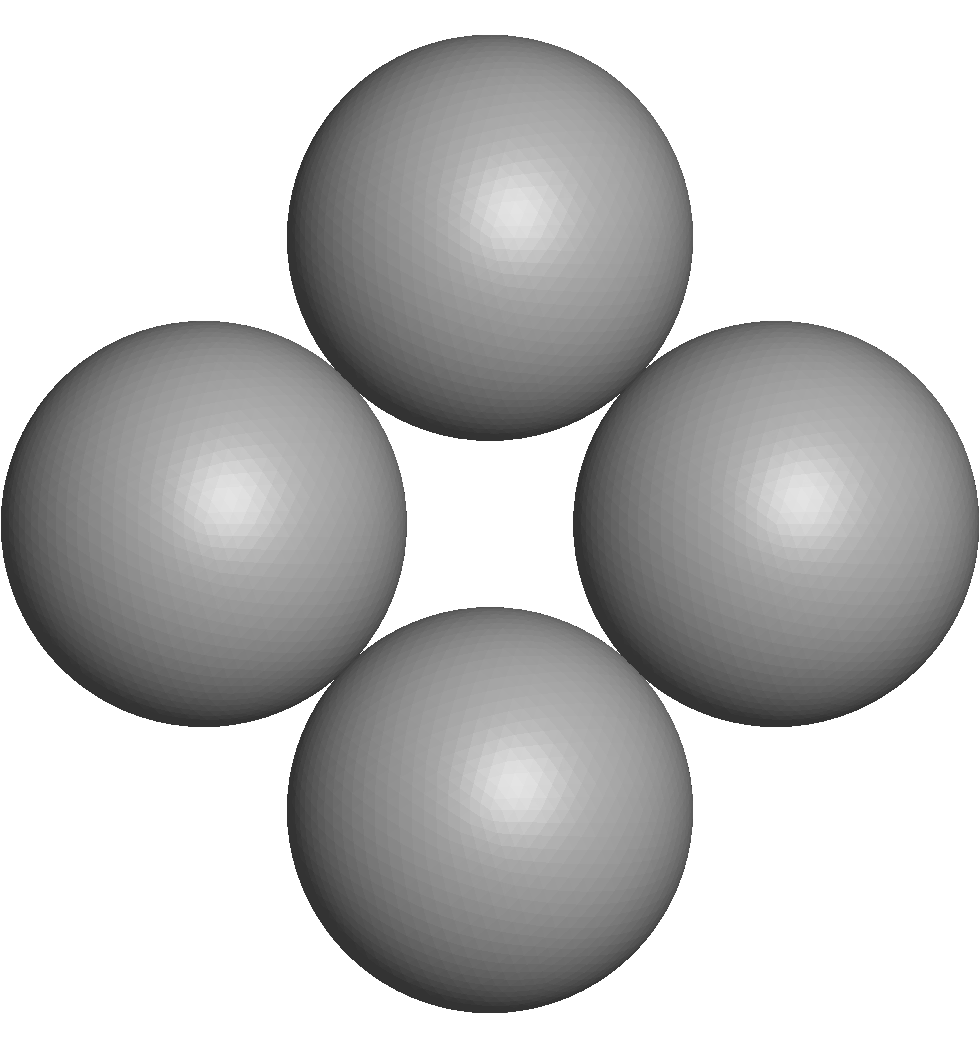} &
    \includegraphics[scale=0.08]{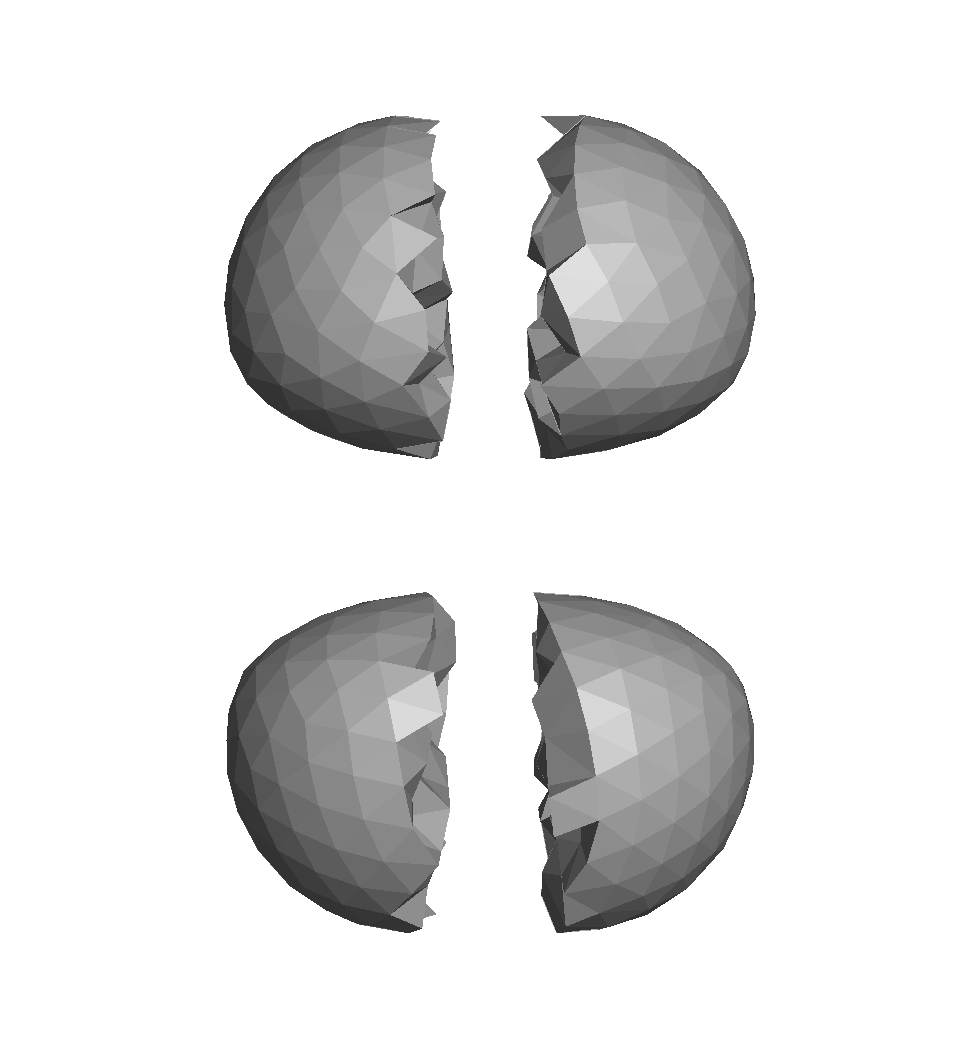}
    \includegraphics[scale=0.08]{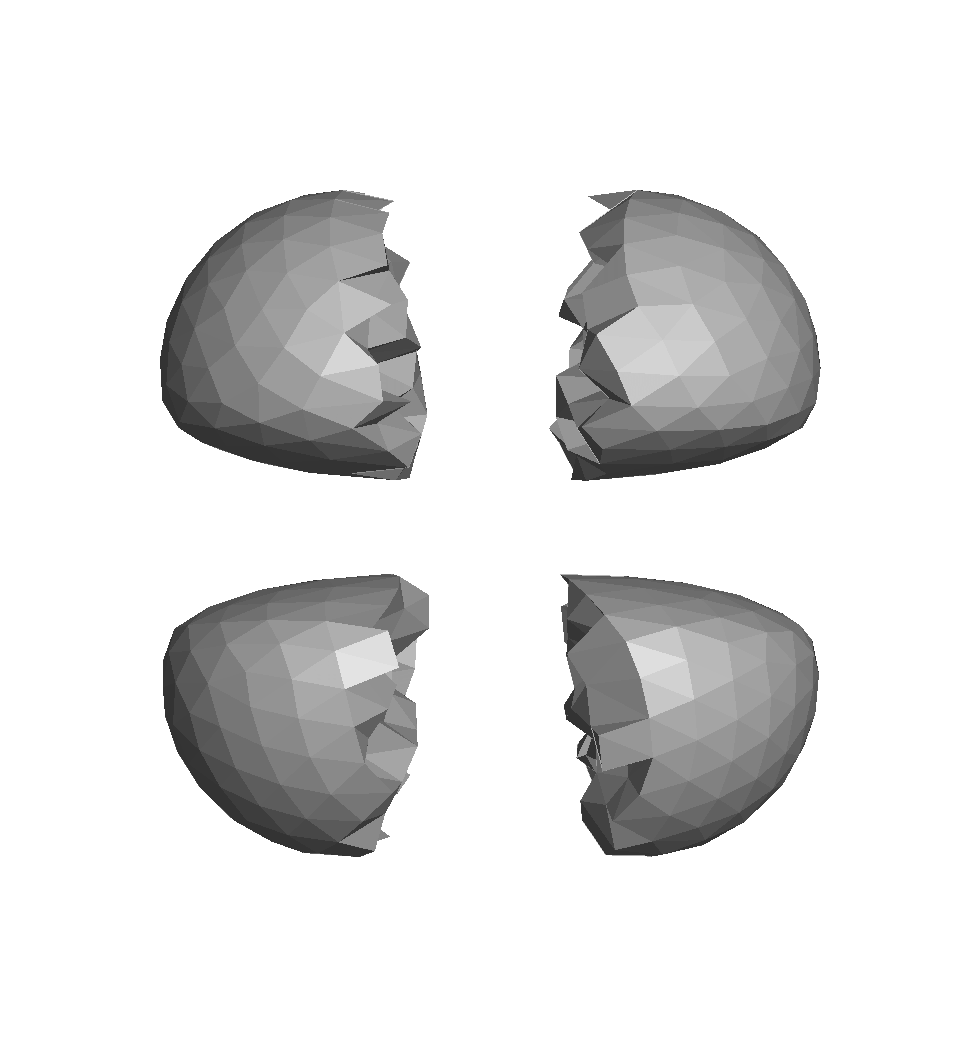}
    \includegraphics[scale=0.08]{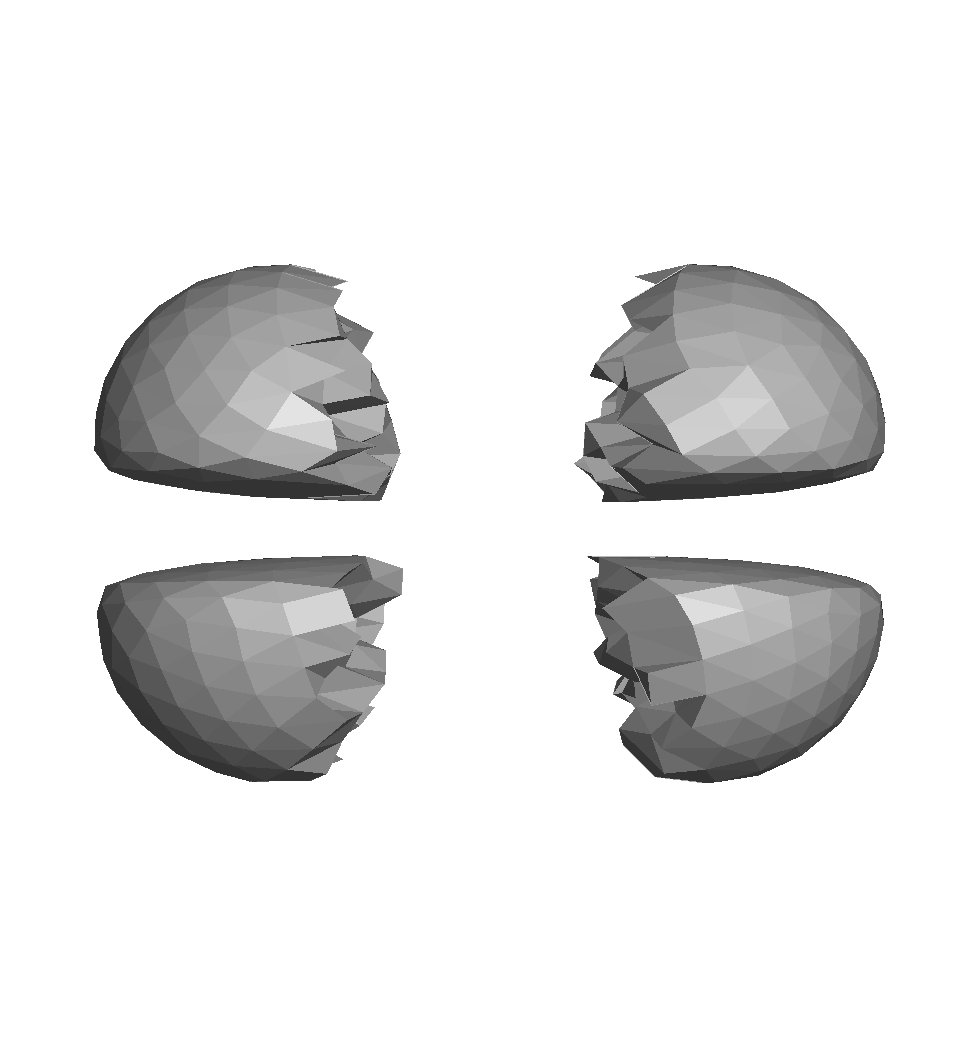} \\
    \hline
  \end{tabular}
  \caption{Interpolation between two horizontally aligned spheres and two
    vertically aligned spheres. Like in 2-d, power cells align with the
    discontinuities in the case of our algorithm, allowing for accurate tearings,
    while with Levy's algorithm, in each case at least one of the tearings is
    not accurately captured.}\label{fig:3d-shape2}
\end{figure}

Figure~\ref{fig:3d-shape2} illustrates the importance of symmetry in 3-d as well:
one can observe the cells aligning with the discontinuities on our algorithm's
output, resulting in two rather accurate tearings, while the interpolations from
Lévy's algorithm always display a ragged tear in at least one of the
discontinuities, in a similar manner as what could be observed on
Figure~\ref{fig:shape9}.

\begin{figure}[!h]
  \centering
  \begin{tabular}{|c|c|c|}
    \hline
     & Input measures & Interpolation \\
    \hline
    \rotatebox[origin=l]{90}{\scriptsize{Our interpolation} } &         \includegraphics[scale=0.08]{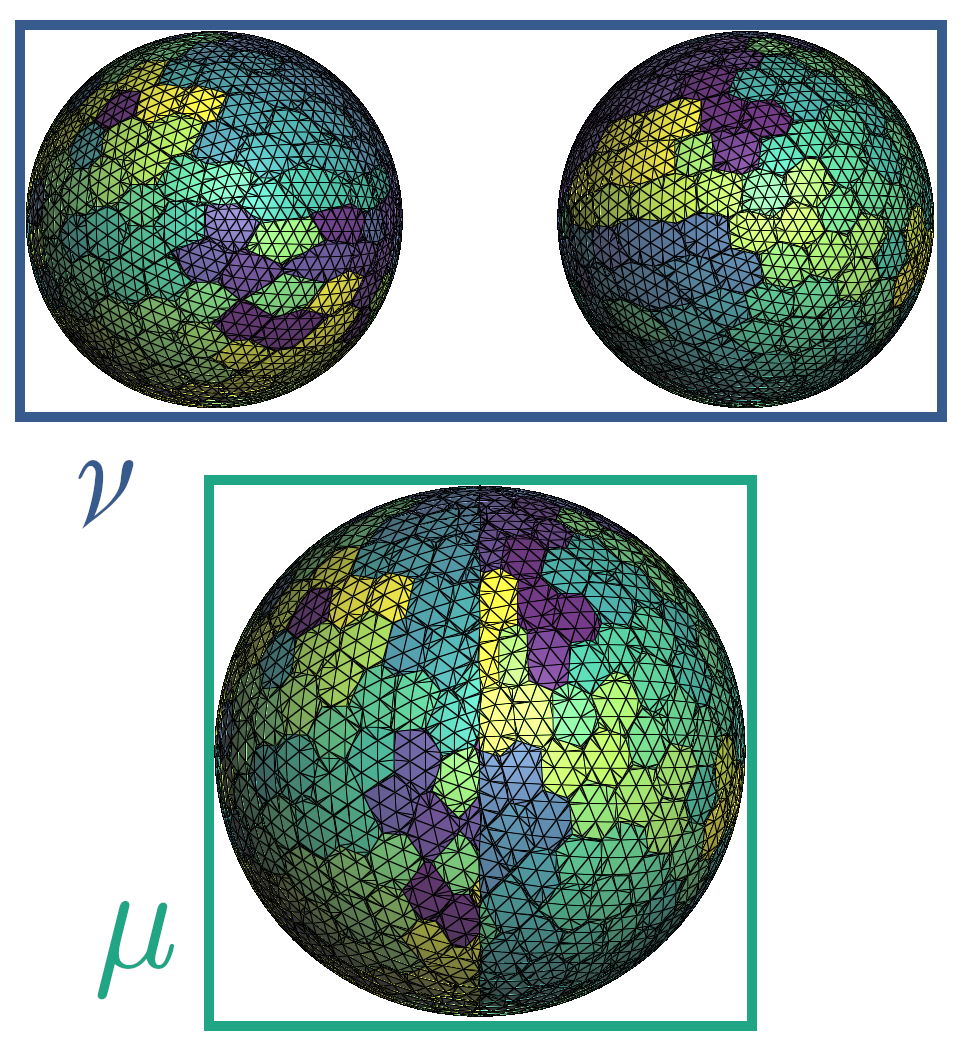} &
    \includegraphics[scale=0.08]{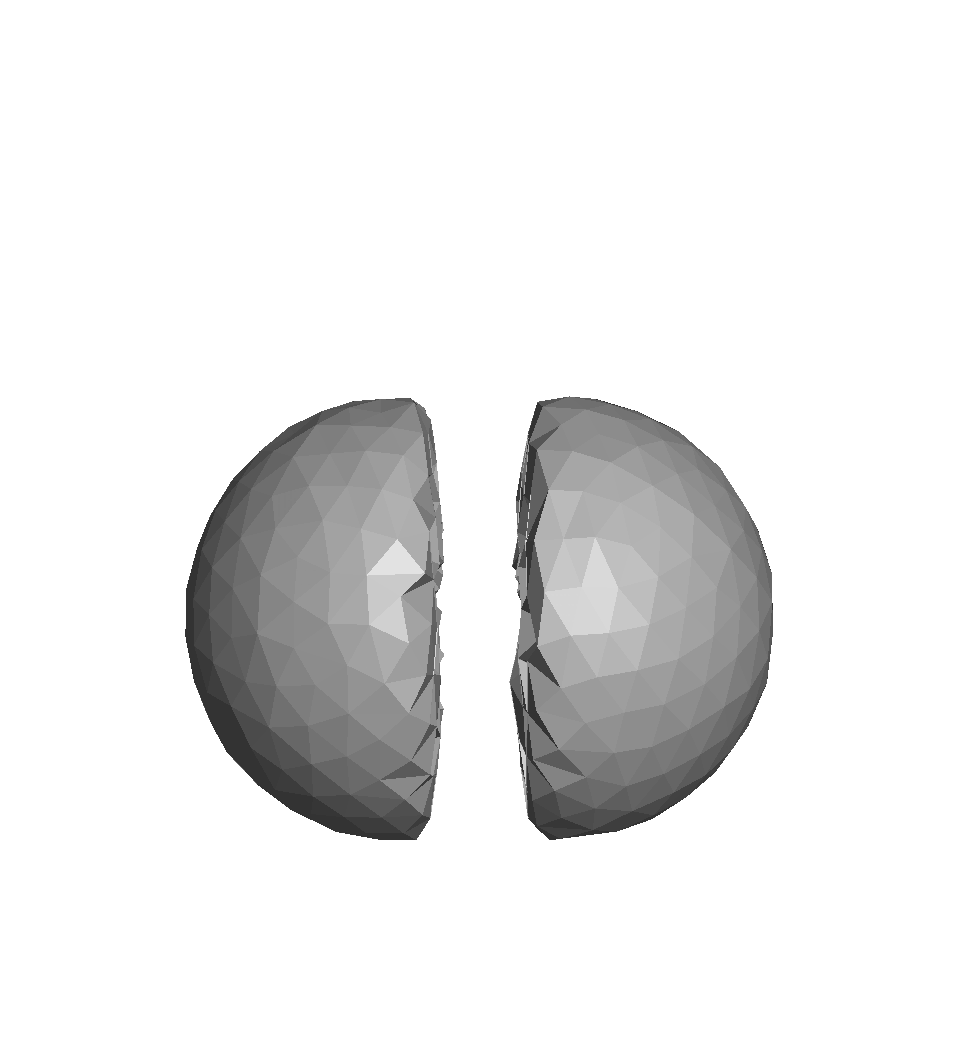}
    \includegraphics[scale=0.08]{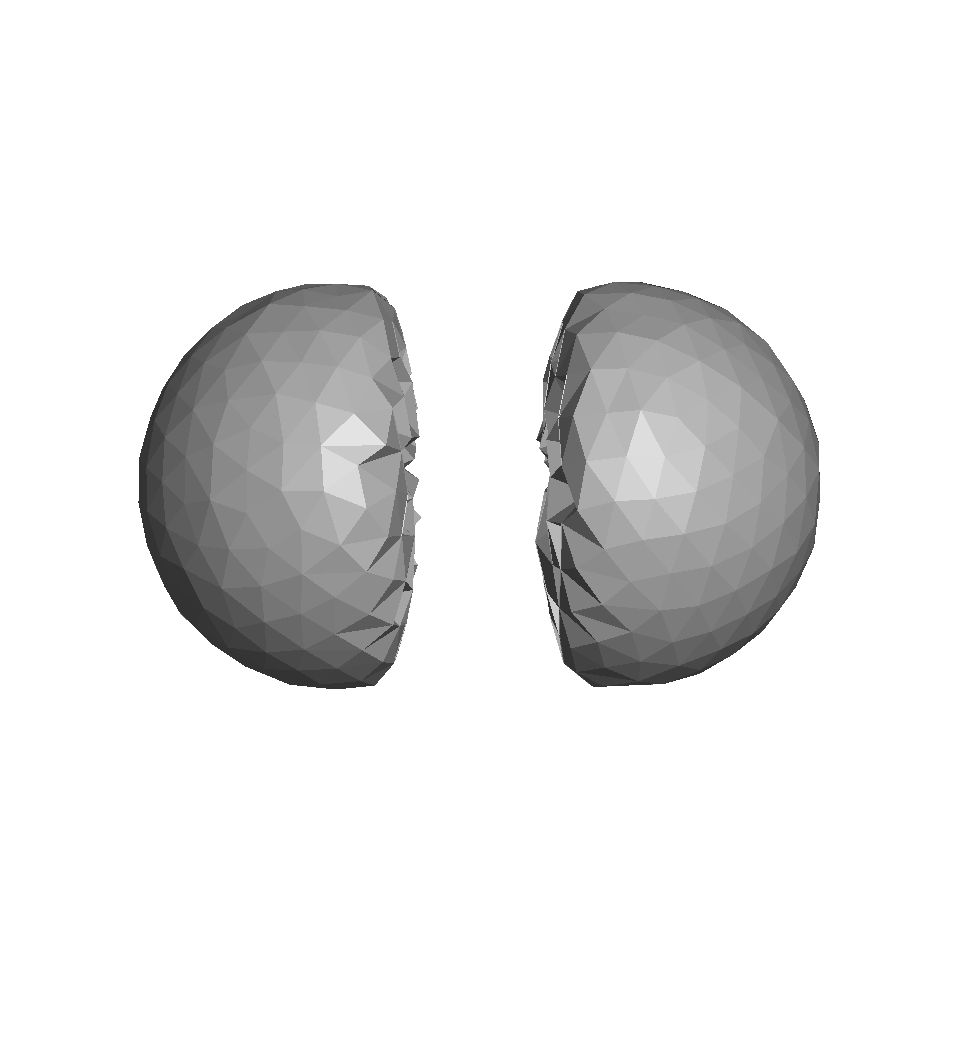}
    \includegraphics[scale=0.08]{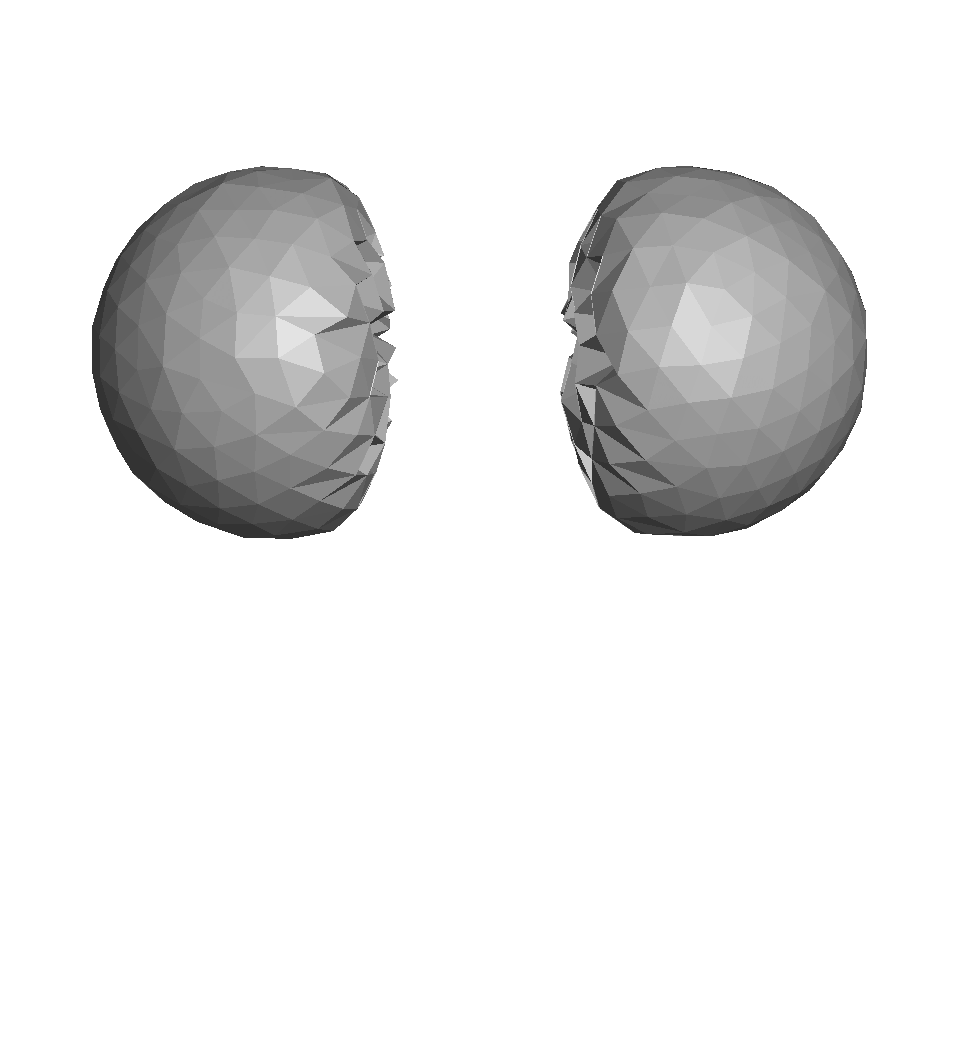} \\
    \hline
    \rotatebox[origin=l]{90}{\scriptsize{semi-discrete~\cite{bl-semi-discrete}}}
    \rotatebox[origin=l]{90}{\scriptsize{$\mu$ continuous}}
    \rotatebox[origin=l]{90}{\scriptsize{$\nu$ sampled}} &
    \includegraphics[scale=0.08]{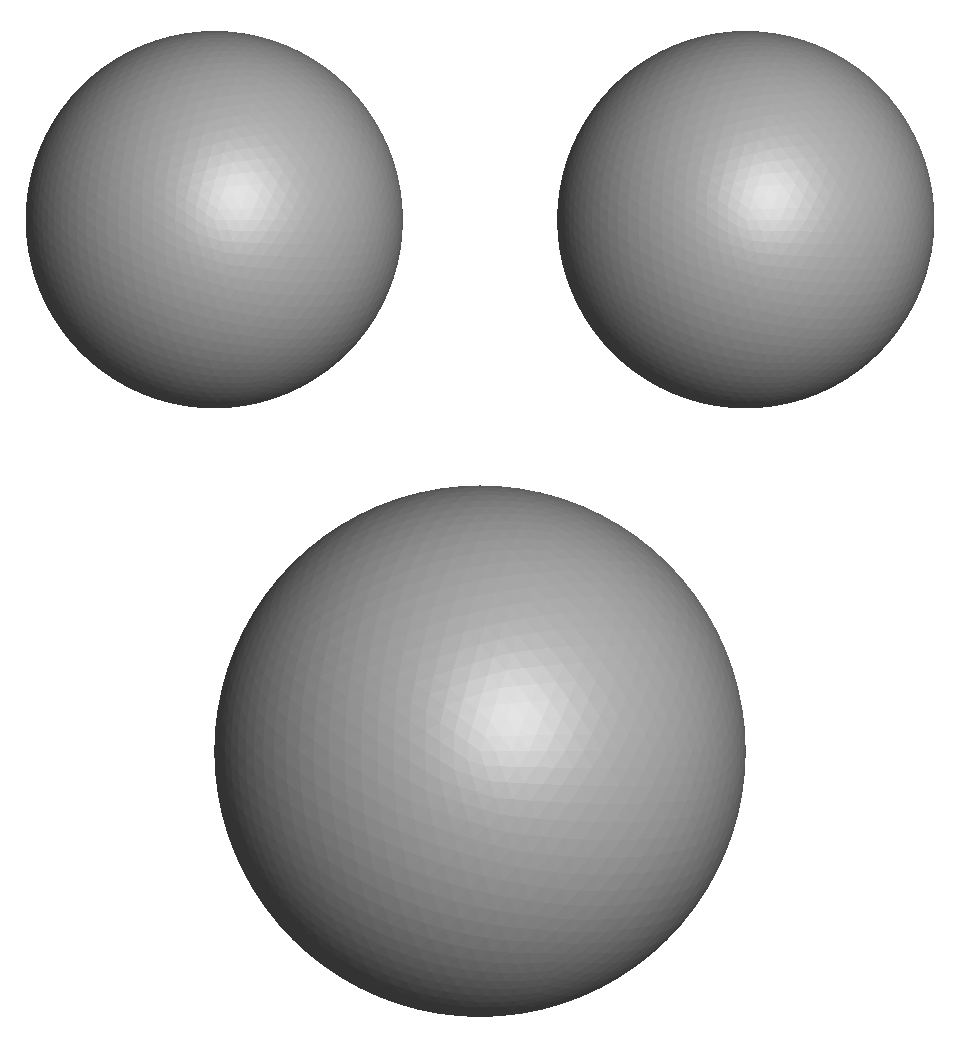} &
    \includegraphics[scale=0.08]{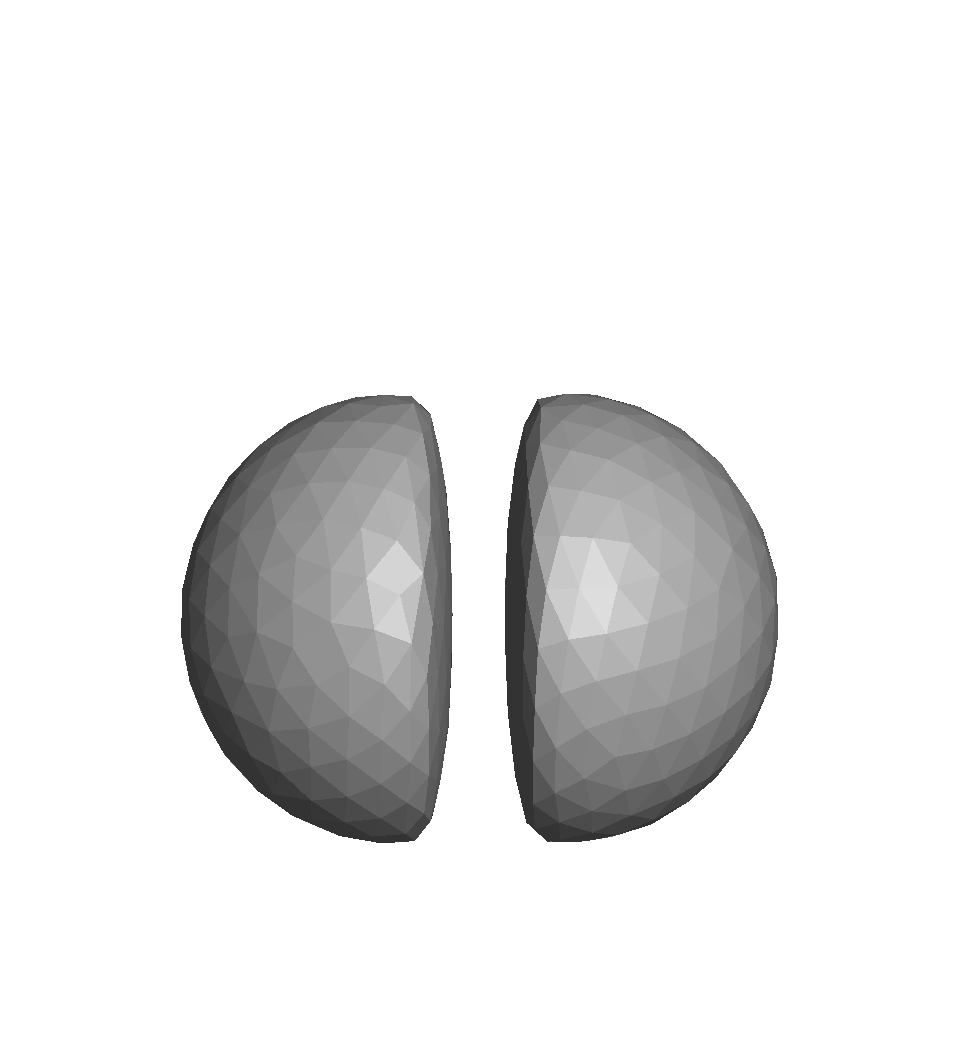}
    \includegraphics[scale=0.08]{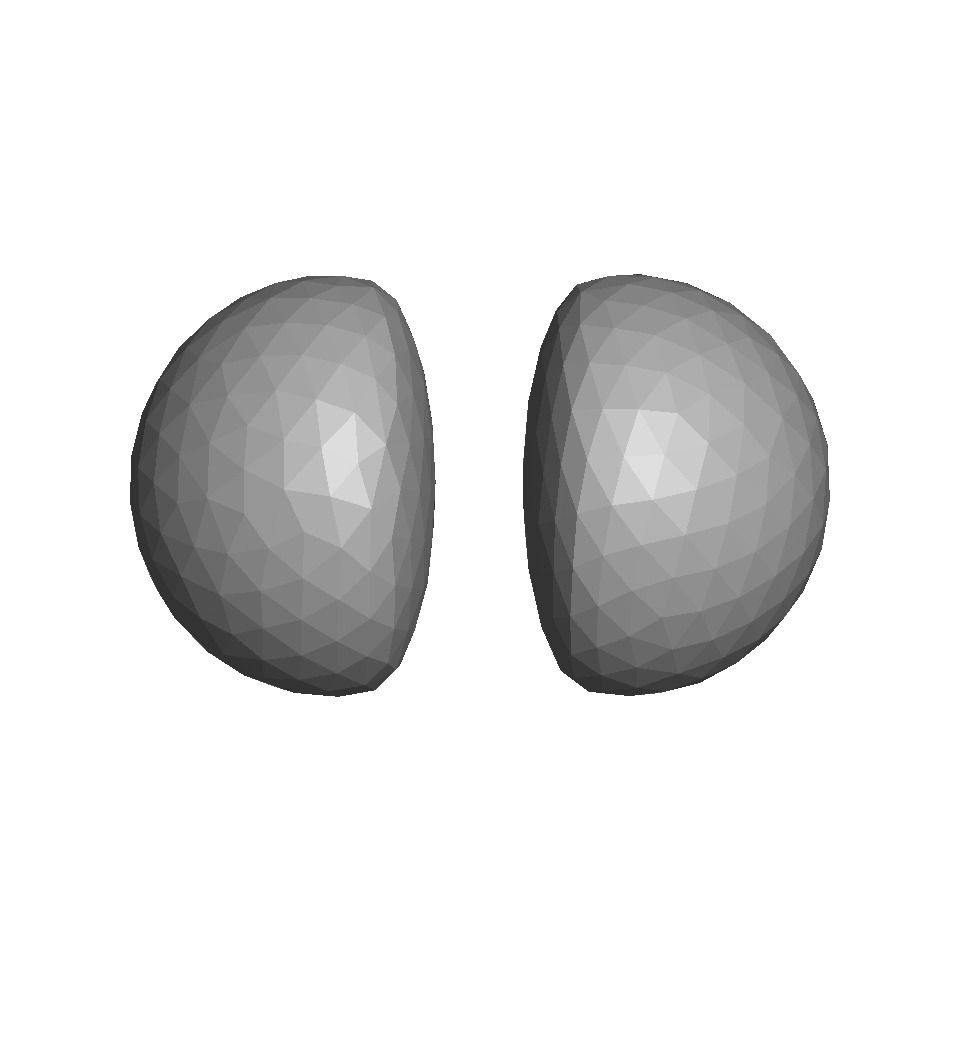}
    \includegraphics[scale=0.08]{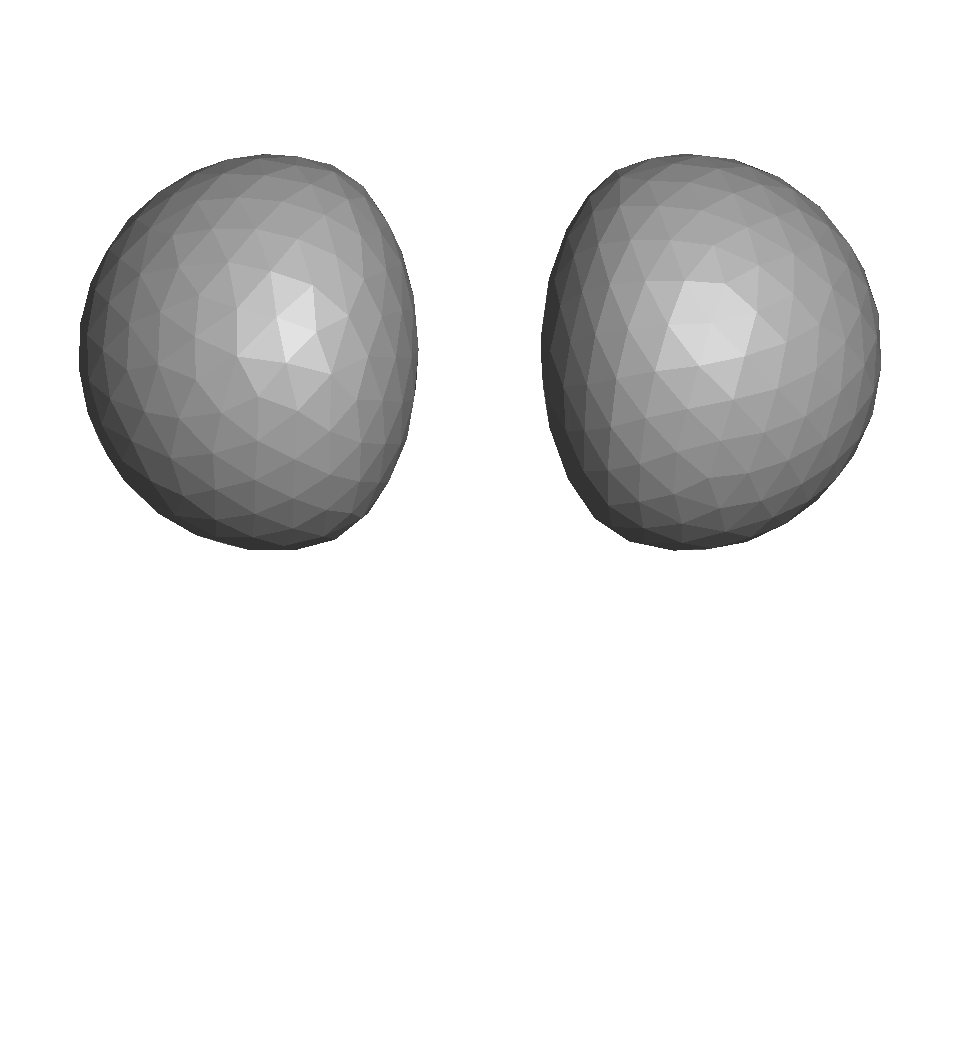} \\
    \hline
     \rotatebox[origin=l]{90}{\scriptsize{semi-discrete~\cite{bl-semi-discrete}}}
     \rotatebox[origin=l]{90}{\scriptsize{$\mu$ sampled}}
     \rotatebox[origin=l]{90}{\scriptsize{$\nu$ continuous}} &    
    \includegraphics[scale=0.08]{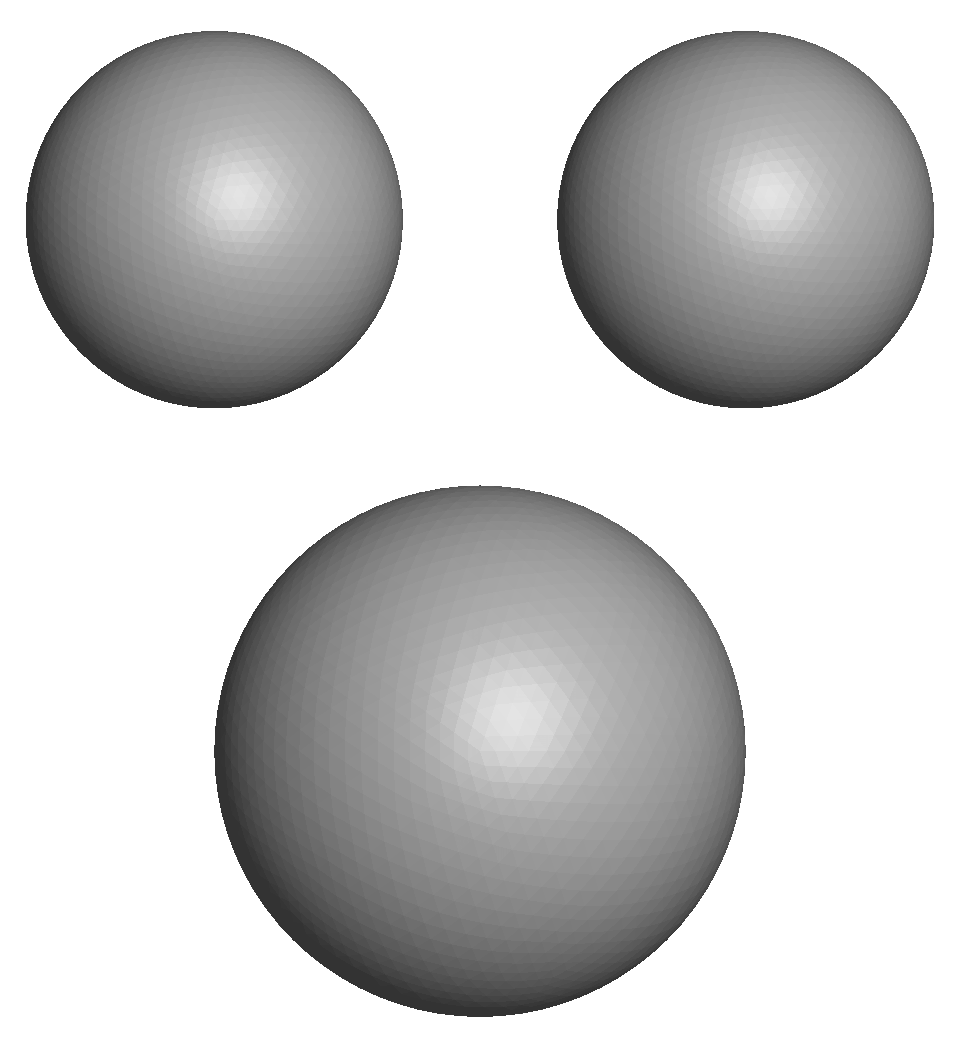} &
    \includegraphics[scale=0.08]{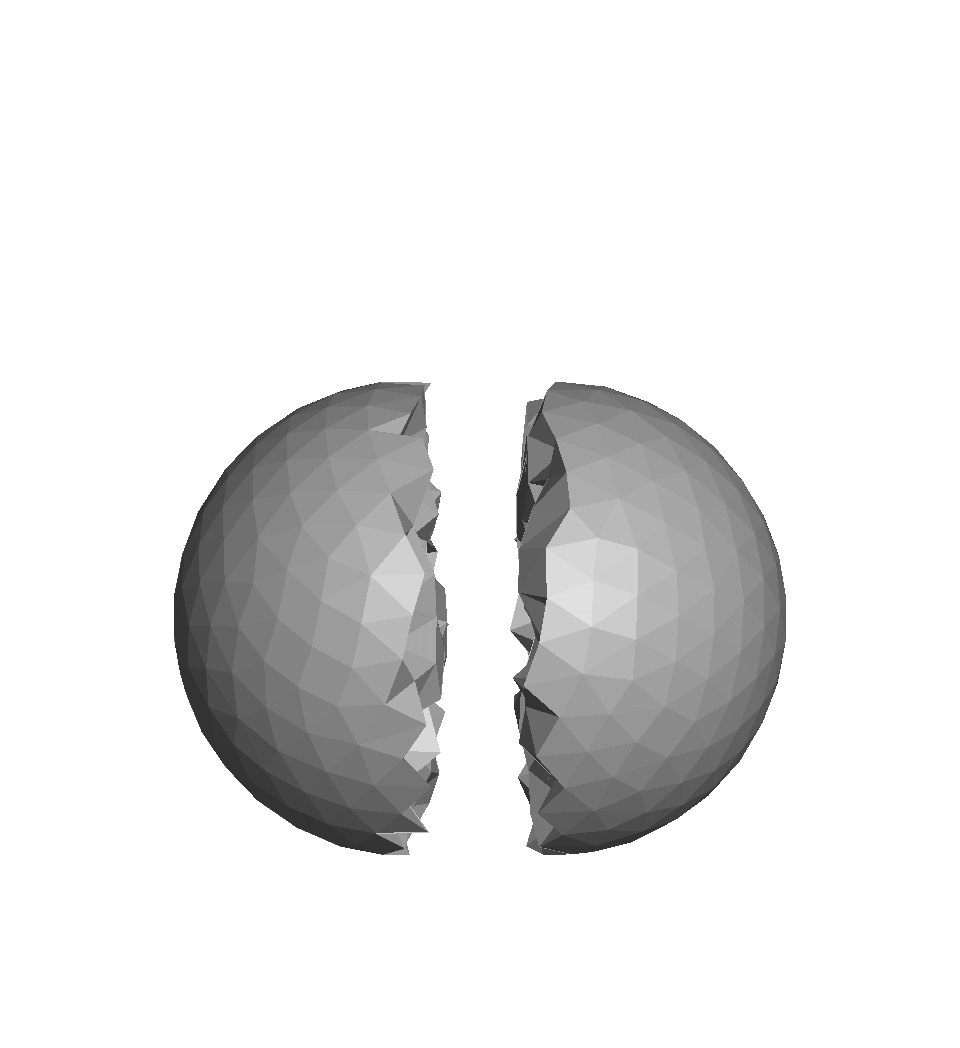}
    \includegraphics[scale=0.08]{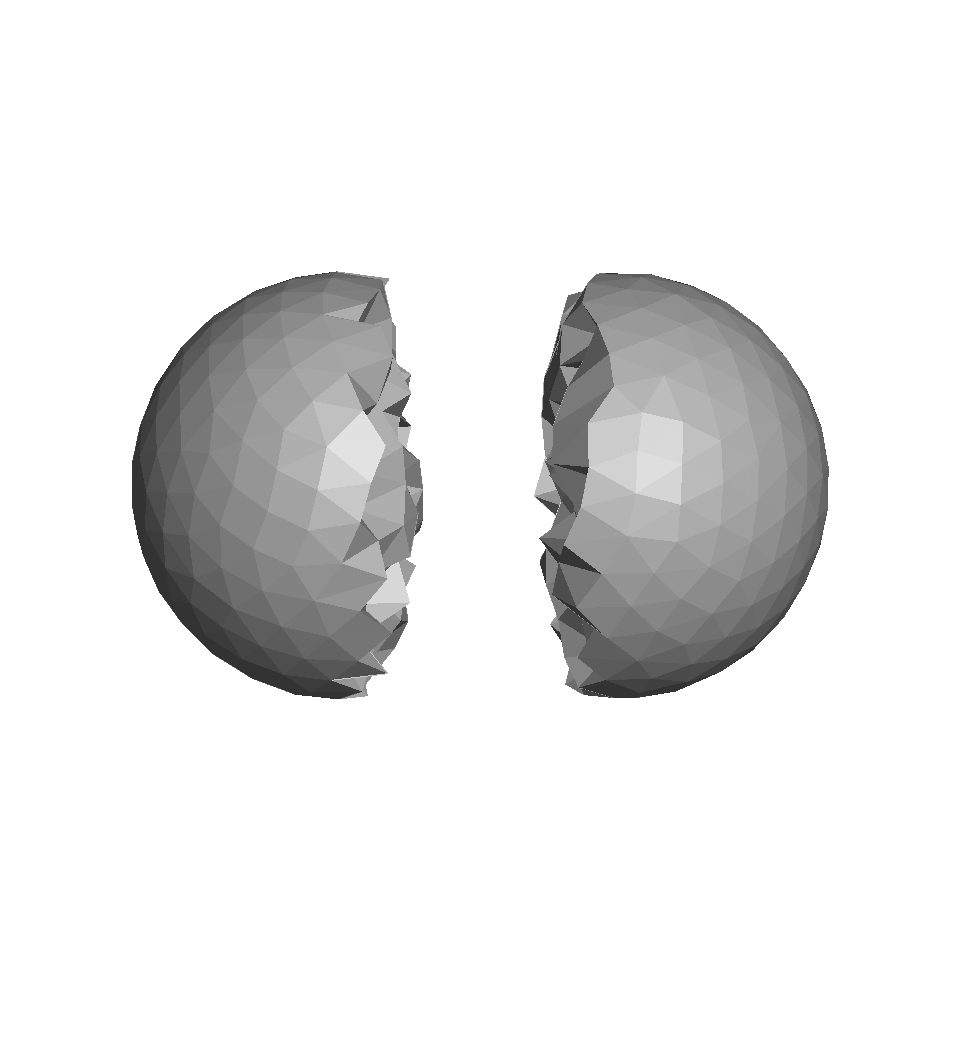}
    \includegraphics[scale=0.08]{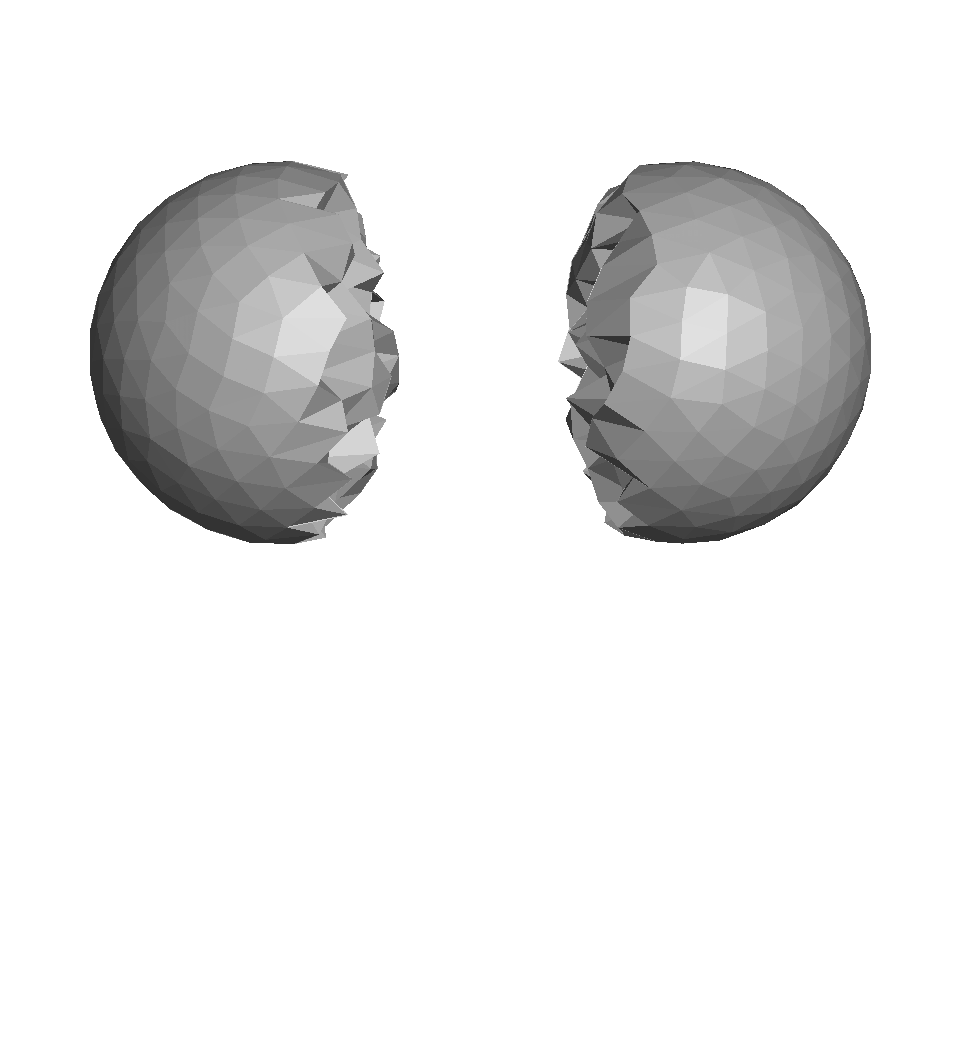} \\
    \hline
  \end{tabular}
  \caption{Interpolation between a single sphere and two smaller spheres. The
    power cells align with the discontinuity with our algorithm but only do
    so in one case with Levy's algorithm.}\label{fig:3d-shape1}
\end{figure}

Figure~\ref{fig:3d-shape1} illustrates a similar phenomenon to the one observed
on~\ref{fig:3d-shape2}: the power cells of the transport map align with the
discontinuity on our algorithm's output, while it is only the case with Levy's
algorithm when the discontinuity is on the target measure.

\subsection{Quantitative evalutation}

In order to evaluate the accuracy of our interpolations, we compute a ground
truth interpolation using Lévy's algorithm with $100k$ samples and compare the
Hausdorff distances between our algorithm with $10k$ samples and
the ground truth, and the Hausdorff distance between Lévy's algorithm with $10k$ samples and the ground truth.
We report the results in Table~\ref{fig:hausdorff}.
Our rationale for using Lévy's algorithm with a dense sampling as a ground truth is that the tearing effects are less critical when the sampling density is high. It might favor in theory Lévy's algorithm, but our quantitative results show that even with this small bias, our method outperforms Lévy's algorithm, when the sampling density is lower.
Indeed, we observe that our results are systematically more precise than Lévy's, which corroborates the visual observations.

\begin{table}[!h]
  \centering
  \begin{tabular}{|l|l||c|c|c|c|c|}
    \hline
    Shape & Algorithm & 0 & 0.25 & 0.5 & 0.75 & 1 \\ \hline
    Six-pointed star to & Levy's & 0.078 & 0.076 & 0.074 & 0.076 & 0.077\\
    six-pointed star~\ref{fig:shape7} & Ours & 0.028 & 0.033 & 0.038 & 0.044 & 0.045 \\ \hline
    Two disks to & Levy's & 0.025 & 0.033 & 0.025 & 0.021 & 0.016\\
    two disks~\ref{fig:shape9} & Ours & 0.0086 & 0.017 & 0.0098 & 0.0094 & 0.0072 \\ \hline
    One disk to & Levy's & 0.022 & 0.022 & 0.025 & 0.026 & 0.016\\
    three disks~\ref{fig:shape2} & Ours & 0.0098 & 0.0089 & 0.0017 & 0.0092 & 0.0064 \\ \hline
  \end{tabular}
  \caption{Hausdorff distances between meshes the interpolated meshes at different
    interpolation steps, }\label{fig:hausdorff}
\end{table}

\section{Limitations and discussion}

\begin{figure}[!h]
    \centering
    \includesvg[width=0.4\textwidth]{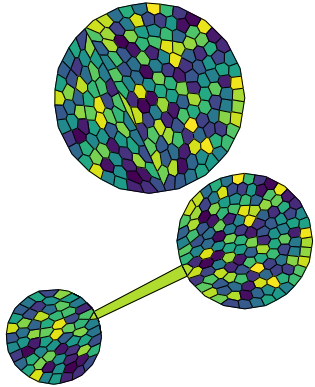}
    \caption{Overlapping cell in converged transport maps}
    \label{fig:overlapping}
\end{figure}

It can sometimes happen that, even when the algorithm has seemingly converged, one cell 
overlaps between connected components, as shown in figure~\ref{fig:overlapping}.

This is usually due to the ratio between the different connected components' areas being
incompatible with the number of cells: for example, a domain composed of two connected 
components of same masses will have to split a cell between both components in order to 
fulfill the mass constraints of optimal transport. This calls for a wise choice of the 
exact number of samples used. 

A limitation of our algorithm lies in its computation time.
To alleviate this, a possibility is to use Newton's method as in Kitagawa, Mérigot and Thibert~\cite{kmt}
to compute semi-discrete optimal transport maps.
However, this raises initialization concerns: the weights of the power diagram have to be initialized in such a
way that there are no empty cells. As a consequence, whenever we deal with two
non geometrically identical domains (which constitutes the vast majority of our
use cases), we cannot initialize our power diagrams as Voronoi diagrams. This
forces us to use some kind of initialization procedure, such as the ones
described in~\cite{meyron-init}, but they have been unsatisfactory in practice.
Improving computation times remains an open question.

\section{Conclusion}

We introduced in this article a novel approach of using semi-discrete optimal transport to
approximate displacement interpolation, by coupling two semi-discrete transport maps 
through the barycenters of their cells. We presented a fixed-point algorithm, following a 
classical alternating pattern, to compute such coupled transport maps, and showed how it
empirically converged to the specified setting. We finally showed how the coupled 
transport maps present strong geometric similarities that allow us to construct an
accurate approximation of displacement interpolation by simply linearly interpolating 
between their vertices. We also observe that the cells of the transport maps accurately 
align with topological discontinuities in the measures' supports, and that this property
is transmitted to the interpolation.

\bibliographystyle{siamplain}

\end{document}